# The Force Singularity for Partially Immersed Parallel Plates

Rajat Bhatnagar and Robert Finn

**Summary: In earlier work, we provided a general description of the forces of attraction and repulsion, encountered by two parallel vertical plates of infinite extent and of possibly differing materials, when partially immersed in an infinite liquid bath and subject to surface tension forces. In the present study, we examine some unusual details of the exotic behavior that can occur at the singular configuration separating infinite rise from infinite descent of the fluid between the plates, as the plates approach each other. In connection with this singular behavior, we present also some new estimates on meniscus height details.**

This paper continues earlier investigations by its authors and by others, into a problem going back to the (unsuccessful) 17th Century attempts by Edme Mariotte to characterize horizontal force interactions among floating bodies. Over a century later Laplace [7] considered the phenomenon in the context of then newly established surface tension theory, and observed that essential simplifications appear for the particular configuration of two infinite parallel vertical plates, partially dipped into an infinite liquid bath. The present work, along with its predecessors [1,2,3], continues in the spirit introduced by that author. Figure 1 (taken from [1]) illustrates the general reference configuration. In common with [1,2,3], our study goes beyond [7] in varying respects, prominently in its emphasis on structure of the exotic singular behavior at the configuration in which the contact angles on the sides of the plates facing each other are supplements $\gamma_1 + \gamma_2 = \pi$, with plate separation tending to zero. For reasons not clear to us, Laplace ignored that singularity in his writings, although he did remark the configuration as having special interest. In the present study we focus attention on details of that singular behavior, and seek initial clarification of its striking and exotic structure.

The material of this paper relates closely to that of Aspley, He and McCuan [4], who point out that for any fixed $\gamma_1$, $\gamma_2$ the curves of force vs. plate separation must be one of three clearly delineated types. In the "Addenda" paper [10] just following the present one, we seek to relate the differing kinds of information provided in the two approaches. The present work continues in our earlier perspective, relating fluid height profiles to plate separation and using selected "barrier" curves as guideposts to distinguish behavior patterns for varying regimes of contact angle pairs $\gamma_1$, $\gamma_2$. As shown in [1], the angles $\gamma_{11}$, $\gamma_{22}$ on the outer sides of the plates have no influence on the (equal) forces exerted by the plates on each other.

In what follows, we use in essential ways material from our earlier papers [1,2,3], and we assume some familiarity with the notation and central results of those references.

The questions to be addressed derive from classical literature and from observations established in [1,2,3,9]: we know from these citations that the surface profile exists as a formal solution of the basic *capillarity equation*,

$$\operatorname{div}\left(\frac{\nabla u}{\sqrt{1+|\nabla u|^2}}\right) = \kappa u, \quad \kappa \equiv const$$

and from [5] that it is the same in every vertical plane orthogonal to the two plates; thus we may restrict attention to a representative such plane, which we designate as the (*x,u*) plane, with the profile described by *u = f(x)*. By [9] the surface profile between the plates meets them in "contact angles" $\gamma_1$, $\gamma_2$, depending only on the materials and in no other way on the conditions of the problem, see Figure 1.

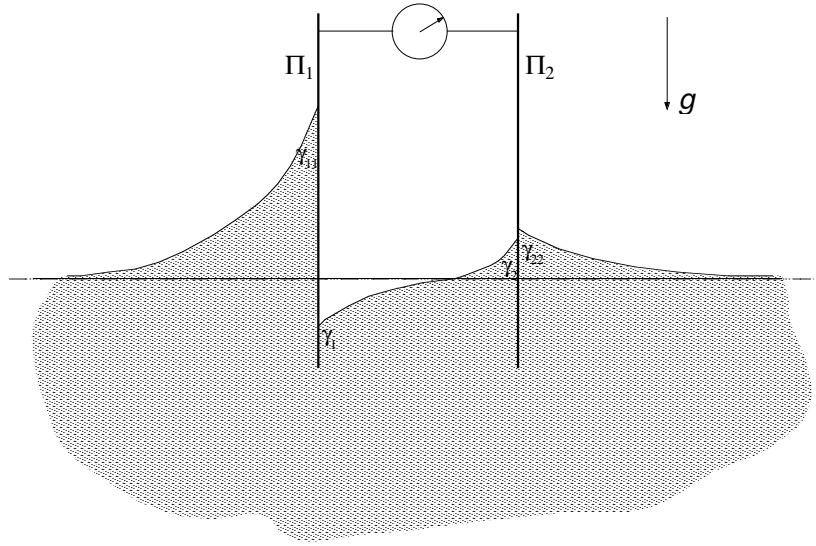

**Figure 1. The underlying configuration. Each plate may be of different materials on its two sides, thus leading to four distinct contact angles. Depending on circumstances the plates may attract or repel each other.**

If $\gamma_1$, $\gamma_2$ are supplementary $\gamma_1 + \gamma_2 = \pi$, then the profile between the plates will be symmetric relative to reflection in a surface point located on the *x* – axis midway between them, and we can show that the net horizontal force density between the plates will be repelling for every plate separation 2*a* > 0. Denoting by ***F*** the force density per unit length (taken as negative for repelling forces) on the contact lines with the plates, ***F*** will decrease monotonically as *a* decreases, and there exists

$$\boldsymbol{F}_0 = \lim_{a\to 0}\boldsymbol{F} = -2\sigma(1-\sin\gamma_1) = -2\sigma(1-\sin\gamma_2) \qquad (1)$$

This limit provides in fact an upper bound on the force magnitude for all achievable symmetric configurations with the same γ values; that bound is not achieved by any surface in the class considered, but is nevertheless approached (from below) as the plates approach each other. (For a more general statement see Lemma 2 below).

*The indicated behavior is remarkable in that for every configuration with $\gamma_1 + \gamma_2 \neq \pi$, there holds $\mathbf{F}_0 = \lim_{a \to 0} \mathbf{F} = +\infty$, the positive sign denoting an attracting force.* We devote the remainder of this paper to a more complete description of this discrepant behavior, which held for us some surprises.

## DISCUSSION

We start with a symmetric configuration arising from supplementary contact angles (**III** in Figure 2) and with general initial plate separation $2a$, and introduce two ranges of comparison configurations, obtained by varying the contact angle $\gamma_1$. We will then allow $2a$ to decrease, and compare the limiting behaviors of the neighboring configurations with that of the symmetric ones. As in our previous work, we make an initial choice $0 \leq \gamma_2 < \pi/2$, which remains fixed throughout our discussion. The remaining case $\pi/2 < \gamma_2 \leq \pi$ reduces to that one via formal reflection in the $x$ – axis (under which the governing equations (2$a,b$) are invariant), and is thus implicitly covered by the present discussion.

Figure 2 (taken essentially from [3]) offers sketches of "*barrier solutions*" in three different configurations, corresponding to three qualitatively distinct choices of plate separation. These barriers control the qualitative behavior of all solutions of the basic equation (2$a,b$) for the free surface height $u(x)$ between the plates, corresponding to prescribed contact angles $\gamma_1, \gamma_2$ with the plates. All curves of Figure 2 are solutions arising from a common angle $\gamma_2$, with $0 \leq \gamma_2 < \pi/2$.

Figure 2$a$ displays the "*wide*" plate separation. All curves represent solutions of the first order system

$$(\sin \psi)_x = \kappa u$$
$$u_x = \tan \psi$$
(2$a,b$)

Here $\psi$ is the inclination of the solution curve $u(x)$, so that the left side of (2$a$) becomes the (signed) planar curvature of that curve; $\kappa$ is the *capillarity constant* $\kappa = \rho g/\sigma$: $\rho$ = (positive) density difference, $g$ = gravitational acceleration, $\sigma$ = surface tension of the fluid/fluid interface. The symmetric solution **III** meets $\Pi_1$ in angle $\pi - \gamma_2$. *The "wide" plate separation is defined as any configuration for which the solution $\mathbf{IV_0}$ (which is attached to $\Pi_2$ and independent of plate separation) does not extend to $\Pi_1$, as indicated in the figure.* Note that **IV** is defined by its contact angles $\gamma_2$ with $\Pi_2$ and $\gamma_1 = \pi$ with $\Pi_1$, and does depend strongly on plate separation.

We allow the plates to come together by keeping $\Pi_2$ fixed and moving $\Pi_1$ toward it. In doing so, the curve **IV** must move downward in order to make its curvatures more negative, so as to move its vertical point toward $\Pi_2$ and thus maintain its prescribed tangential contact with $\Pi_1$; as the plates approach each other, the entire curve **IV** must move below any prescribed height. The curve $\mathbf{IV_0}$ however does not change; at the separation $2a_0$ corresponding to the width of $\mathbf{IV_0}$, those two curves become identical. For any smaller separation, **IV** lies entirely below **IV₀**. Corresponding to the further interval in which **IV** still lies above **V** we obtain the "*intermediate*" configuration of Figure 2*b*.

By decreasing the separation still further, **IV** at first coincides with the curve **V** (which is rigidly attached to $\Pi_2$) at a separation $2a_{00} < 2a_0$, and then descends below **V**. We designate the configurations thus obtained as the "*narrow*" configurations. In Appendices 1 and 2 we provide explicit expressions for determining the crossover values $a_0$ and $a_{00}$. (In our earlier paper [3], $a_{00}$ is designated $a_{cr}$.)

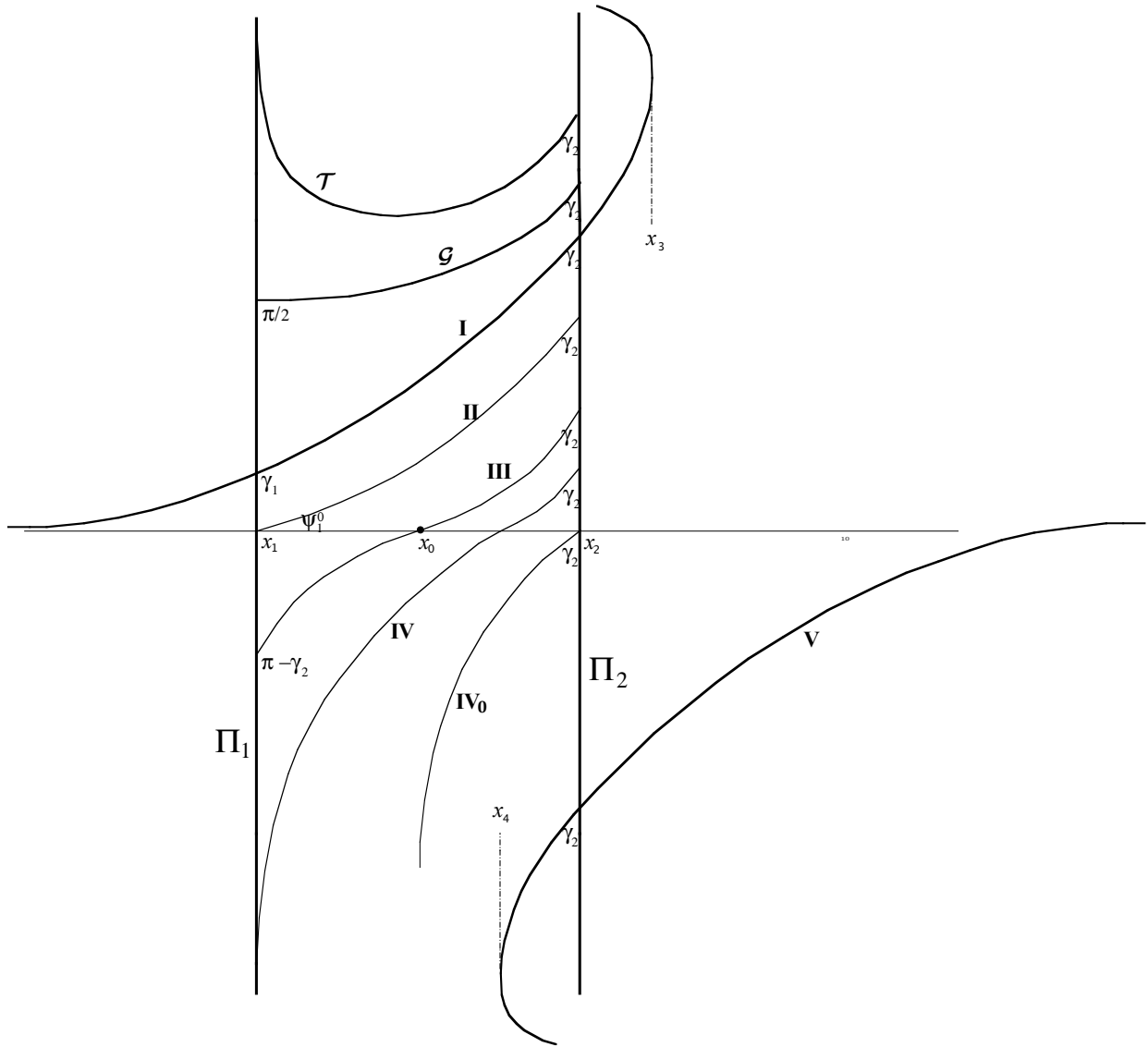

**Figure 2*a*. Sketch of a typical "wide" configuration. IV lies above IV₀.**

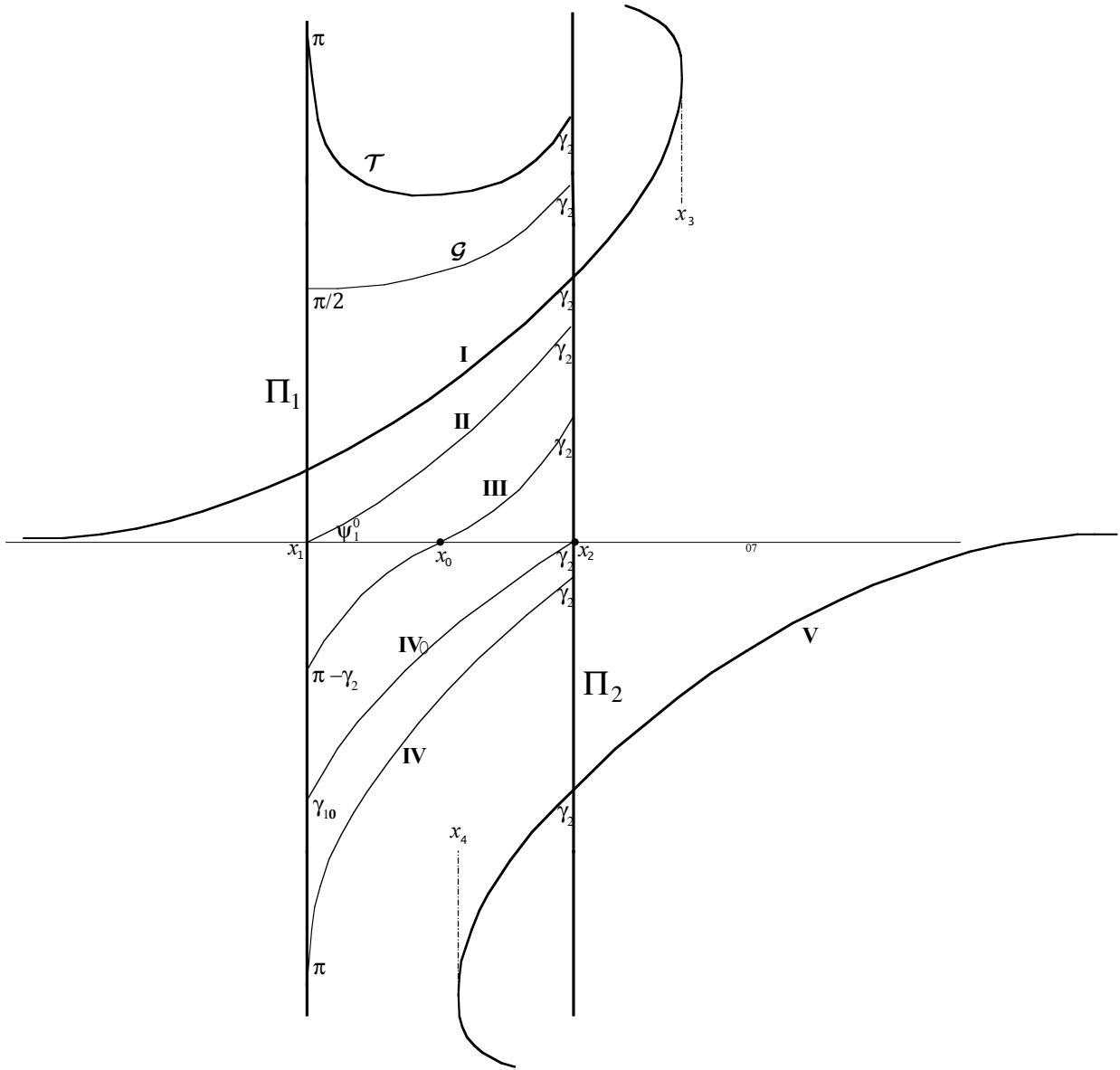

**Figure 2*b*.  Sketch of a typical "intermediate" configuration. IV lies between IV$_0$ and V.**

**Figure 2*c*. Sketch of a typical "narrow" configuration. IV lies below V.**

For later reference, throughout Figure 2 the regions between any two barriers $\Phi$ and $\Psi$ joining both plates will be denoted by symbols of the form $\mathcal{R}_{\Phi-\Psi}$. For example $\mathcal{R}_{\text{I\_IV}}$ denotes the region between the plates cut off by the two barriers **I** and **IV**.

In what follows we introduce the normalized (non-dimensional) coordinates

$$\xi = x/a, \ U = u/a, \ \mathrm{B} = \kappa a^2 ; \tag{3}$$

this leads to the conceptual simplification of joining together elements of "equivalence classes" of solutions that are identical up to similarity. The plate separation $2a$ is now replaced by the non-dimensional $2\sqrt{\mathrm{B}}$. Equation (2*a,b*) becomes

$$\begin{aligned}(\sin\psi)_\xi &= \mathrm{B}U \\ U_\xi &= \tan\psi\end{aligned} \tag{2*a,b}$$

The substance of the discussion further above does not change, and we proceed to introduce two sets of "*neighboring solutions*" to a given family of symmetric solutions **III**: $\mathcal{W}(\xi;\mathrm{B};\gamma_2)$ determined by B as parameter. We emphasize again that $\gamma_2$ has been chosen in the range $0 < \gamma_2 < \pi/2$, and remains unchanged through the discussion.

Note that the choice of non-dimensional coordinates differs from the choice made in [6], reflecting the differing perspectives of the two papers.

**Neighbors:**

1) Given a symmetric $\mathcal{W}(\xi;\mathrm{B}_1;\gamma_2)$ joining plate $\Pi_1$ with $\Pi_2$, with separation parameter $\mathrm{B}_1$, we introduce an "*upper class neighbor*" $U^+(\xi)$ by the requirement that on $\Pi_1$, $\psi = \psi_1^+$, for some $\psi_1^+$ in the range

$$\psi_1^0 < \psi_1^+ < (\pi/2) - \gamma_2 = \psi_2 \tag{4}$$

and on $\Pi_2$, $\psi = \psi_2$. Here $\psi_1^0$ is the inclination of **II** on $\Pi_1$, see Figure 2*a*. This choice defines the solution $U^+(\xi)$ uniquely for any $\psi_1^+$ satisfying (4), and permits exactly those solutions joining the two plates, which lie between the barriers **II** and **III** and meet $\Pi_2$ in angle $\gamma_2$. Note that this definition depends on plate separation, since the barrier **II** which defines $\psi_1^0$ does so; however for each separation $\mathrm{B}_1 > 0$, the quantity $\psi_2 - \psi_1^0$ is defined and positive, which is all that we will need. $\psi_1^0$ is evaluated in Appendix 3.

Under our hypothesis $0 < \gamma_2 < \pi/2$, the left and right sides of (4) are clearly defined for all plate separations. For the region below **III** on the other hand, account must be taken that for wide separations some barriers do not extend to meet both plates. We proceed as follows:

2) We again start with a symmetric **III**: $\mathcal{W}(\xi;\mathrm{B}_1;\gamma_2)$ joining plate $\Pi_1$ with $\Pi_2$, with separation parameter $\mathrm{B}_1$. We introduce a "*lower class neighbor*" $U^-(\xi)$ by the requirements:

$2_0$) Suppose $B_1 > B_0 = \kappa a_0^2$. Then the configuration will be "wide", and **IV** will lie above **IV₀** on the interval between the plates for which **IV₀** is defined as a graph. On $\Pi_1$ we choose $\psi = \psi_1^-$, with $\psi_1^-$ arbitrary in the range

$$\psi_2 < \psi_1^- \leq \pi/2. \tag{5}$$

Then $U^-(\xi)$ will be uniquely determined by the angles $\psi_1^-$ on $\Pi_1$ and $\psi_2 = \pi/2 - \gamma_2$ on $\Pi_2$, and will lie in the region between **III** and **IV**, or else be identical with **IV**. Hence it will also lie above **IV₀** over the interval in which that curve is a graph.

$2_{00}$) Suppose $B_1 \leq B_0$. Then we are confronted with an intermediate or narrow configuration, and **IV** either coincides with or lies below **IV₀**, which now extends to meet both plates in well-defined angles. It will be essential for us to restrict ourselves to solution curves lying above **IV₀**, and thus we define $U^-(\xi)$ via the angle $\psi = \psi_1^-$ on $\Pi_1$, with

$$\psi_2 < \psi_1^- < \gamma_{10} - \pi/2 = \psi_{10} \tag{6}$$

where $\gamma_{10}$ is the "contact angle" with which **IV₀** intercepts $\Pi_1$, see Figure 2*b*. We provide a determination of this angle in Appendix 4.

Defining $\gamma_{10} = \pi$ when **IV₀** does not extend to $\Pi_1$, then (6) applies to all cases, and we shall use it in that sense.

One should observe that $\psi_1^+ < \psi_1^-$. The choice of upper index was motivated by position of the curve, rather than size of the angle.

Figure 3*a* provides a heuristic sketch illustrating the construction. Figure 3*b* supports the construction with computer calculations for typical cases of interest. An essential feature is the range of freedom available for each of the choices $\psi_1^+$, $\psi_1^-$. We shall demonstrate that some features of that freedom disappear during the course of (and prior to the completion of) the limiting procedure in which the plates are brought together.

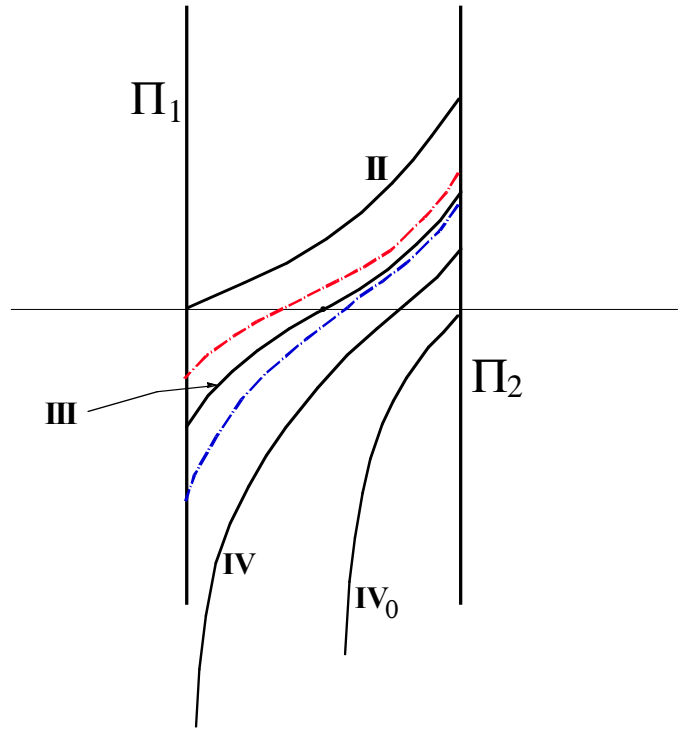

**Figure 3*a*. The ranges of choice for upper class neighbors (red) and for lower class neighbors (blue); illustrated in the case of wide plate separation.**

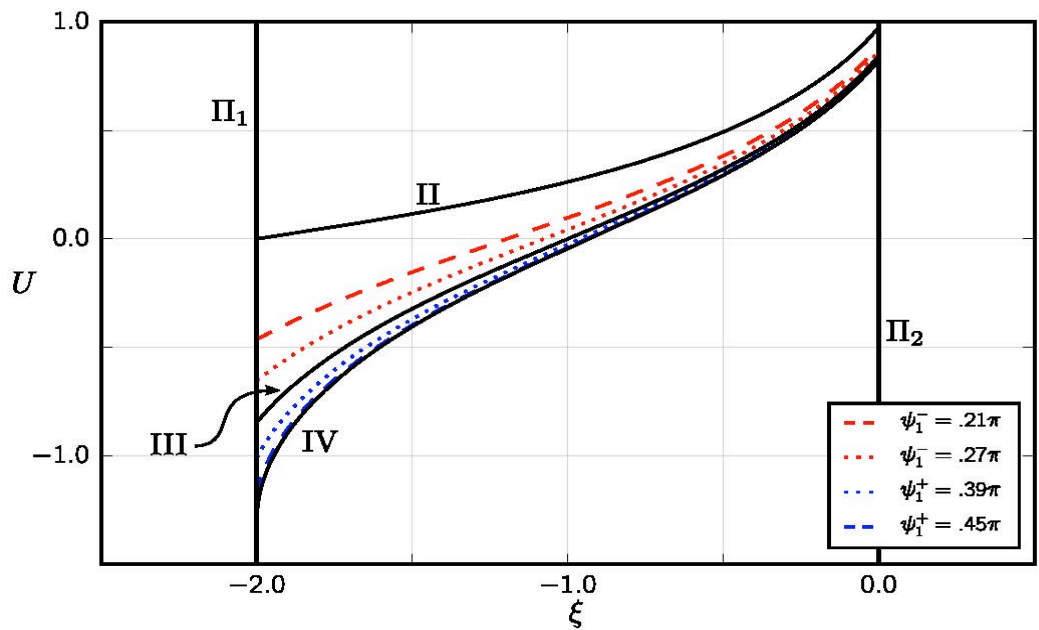

**Figure 3*b*. The configuration of Figure 3*a*; some computer calculations, $\psi_2 = \pi/3$.**

No further restrictions need be introduced. Note that for our idealized purposes the distinctions between "upper class" and "lower class" solutions for fixed $\gamma_2$ and B need not be large. For any given $\mathcal{W}(\xi;B;\gamma_2)$ our conditions in all cases enable us to find "neighbors" from both classes as close to each other as desired, uniformly on the interval between the plates, simply by choosing $\psi_1^-$ and $\psi_1^+$ both close to $(\pi/2) - \gamma_2$. The former solution lies however always strictly below $\mathcal{W}$, the latter one strictly above $\mathcal{W}$; more strikingly, no matter how close to each other the two solutions are for a fixed initial $B$, the asymptotic behaviors as the plates come together can differ greatly in one sense, while remaining completely rigid in another sense, as we proceed to show.

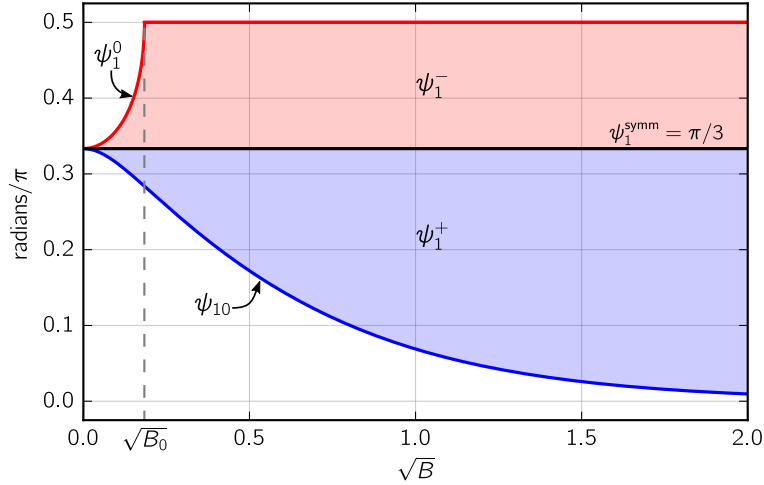

**Figure 4.** Ranges available for $\psi_1^+$ and for $\psi_1^-$, for varying choices of $B$.

Figure 4 illustrates the exact ranges available for $\psi_1^+$ and for $\psi_1^-$, in terms of the parameter B, for the particular choice $\gamma_2 = \pi/6$. The existence of the solutions we have been describing is a special case of the latter author's discussion on the topic in [1]. With regard to our immediate needs, we may state

**Theorem 1.** *As* $\psi_1^+ - \psi_1^- \to 0$ *for any fixed* B, *both the red and the blue solution curves converge uniformly in height to* **III** *over the entire traverse.*

To this effect we use:

**Lemma 1.** *Two distinct solution curves $U(\xi)$ and $V(\xi)$ of (2\*) with the same inclinations at a point p on a common interval of definition cannot intersect on that interval.*

**Proof:** We may suppose $U(p) > V(p)$. Denote the respective inclinations by $\varphi$ and $\omega$. Suppose the curves intersect at an initial such point $q > p$. We obtain from (2*) the contradiction

$$0 < B\int_p^q (U - V)d\xi = \sin\varphi(q) - \sin\omega(q) \leq 0. \tag{7}$$

An analogous proof applies when $q < p$.

Thus neither neighboring curve crosses **III**, since the slopes agree on $\Pi_2$. Since the left side of (2*$a$) is precisely the planar curvature, we see that the height difference between each of the approximating curves and the symmetric curve decreases on moving from $\Pi_1$ to $\Pi_2$; this completes the proof of the theorem.

We now introduce an arbitrary "upper class neighbor" with angles $\psi_1^+$, $\psi_2 = \pi/2 - \gamma_2$, on the two plates. Observe that $\psi_1^+$ is arbitrary in the range $\psi_1^0 < \psi_1^+ < \psi_2$, see Figure 2a. That is, the range of angles available to us is precisely the range of angles on the original barrier **II** at the start of the procedure. Thus there is a point $p_1^+ \in$ **II**, corresponding to a coordinate $\xi_1^+$ between the plates, at which the inclination of **II** is exactly $\psi_1^+$. That means that when the translated $\Pi_1$ has reached the coordinate $\xi_1^+$, the solution will have inclinations at its endpoints exactly those of the restriction to the new and shortened interval, of the curve **II** that appeared in the starting configuration. Thus, denoting that solution by $U(\xi)$ and **II** by $V(\xi)$, we find $\int_{\xi_1^-}^{\xi_2}(U-V)d\xi = 0$. By Lemma 1 the curves do not cross; thus on the interval between $\xi_1^+$ and $\xi_2$, the solution $U$ coincides with **II**.

As a side remark, we note that throughout the described procedure, the solution curve cannot leave the region between **II** and **III**. In fact, the entire curve moves upward away from **III** until the starting point reaches the $\xi$ – axis, cf. the discussion in Sec. 7 of [2]. We do not explicitly need this information for the present context.

A particular consequence is that when the shifted $\Pi_1$ is located at $\xi_1^+$, the entire solution segment joining the two plates is at positive height, and hence its initial point over $\xi_1^+$ has positive height. Since at the start of the procedure the initial point had negative height, we conclude from continuity considerations the existence of an intermediate coordinate $\xi_{10}^+$, at which the initial point of the solution segment lies on the $\xi$–axis. The inclination of that initial point is known to be $\psi_1^+$; recalling from [2] the general force formula for repelling plates

$$\mathcal{F} = -2(1-\cos\psi_0) \tag{8}$$

where $\psi_0$ is the angle made with the $\xi$ axis by a solution curve at the crossing point, and $\mathcal{F}$ is the normalized horizontal force $\mathcal{F} = F/\sigma$, we have achieved:

**Theorem 2.** *As the plate $\Pi_1$ is moved toward $\Pi_2$ with fixed angles $\psi_1^+$ and $\psi_2$, there is an intermediate separation, at which the normalized force $\mathcal{F}^+$ between the plates will be exactly*

$$\mathcal{F}^+ = -2(1-\cos\psi_1^+). \tag{9}$$

We refer in general to a signed force as an *algebraic* force, and note that a negative minimum of an algebraic force is a local maximum for force magnitude.

To Theorem 2 we adjoin:

**Lemma 2.** *If the inclination of a repelling solution curve at a point $(\xi, \eta)$ is $\psi_\eta$, then $|\mathcal{F}| \leq 2(1-\cos\psi_\eta)$, equality holding if and only if $\eta = 0$.*

**Proof:** If $\eta > 0$, then the solution curve is convex upward, from which we conclude that the angle $\psi_0$ of crossing with the $\xi$ axis satisfies $|\psi_0| < |\psi_\eta|$. If $\eta < 0$, then the curve is convex downward, and we obtain the same conclusion. From these two remarks, the assertion follows.

We have shown in [2] that as $\Pi_1$ moves toward $\Pi_2$, the initial point of the solution curve joining the plates, with fixed angles $\psi_1$ and $\psi_2$, moves strictly monotonically upward. Thus Lemma 2 implies

**Theorem 3a.** *The position $\xi_1^{++}$ at which* (9) *holds is uniquely determined*

Also since the inclination $\psi_1^+$ appears (at the starting point) on each of the constructed curves, we find from Lemma 2 and Theorem 3a

**Theorem 3b.** *The position $\xi_1^{++}$ yields uniquely the minimum algebraic force achieved during the procedure.*

Note that this minimum point is *strictly interior* to the interval of the motion. As the plates come still closer together, the algebraic force rapidly increases and tends to positive infinity, no matter how close the initial approximation to **III** may be.

A point to be observed is that the admissible range for $\psi_1^+$ as given by (4) increases with the separation B of the plates. We do not know the extent to which this circumstance reflects reality. With increasing separation, the allowable range for the theorem tends to $0 < \psi_1^+ < \psi_2$. Of course, for values $\psi_1^+$ close to zero the theorem asserts little more than the negativity of the algebraic force.

Note here that by construction, $\psi_1^+ < \psi_2 = \pi/2 - \gamma_2$, which is the value achieved exactly by **III**, and that $|\psi_2 - \psi_1^+|$ can be made as small as desired. Thus, we obtain from Theorems 2 and 3:

**Theorem 4a.** *Let $\psi_{1,j}^+$ be an increasing sequence, satisfying (4) and approaching $\psi_2 = \pi/2 - \gamma_2$. Then the corresponding normalized forces $\mathcal{F}_j$ satisfy $\mathcal{F}_j = -(1 - \cos\psi_{1,j}^+)$, and $\mathcal{F}_j \searrow \mathcal{F}^0 = -(1 - \cos\psi_2)$ as $j \to \infty$. The value $\mathcal{F}_j$ is achieved at a coordinate $\xi_j$, strictly between the plates, with $\xi_j \nearrow \xi_2 =$ the position of $\Pi_2$. The lower bound $\mathcal{F}^0$ is not achieved during this procedure.*

This behavior occurs despite that with further decreasing separation, the force for the neighbor becomes unbounded attracting, and the force for **III** remains repelling and bounded. By Lemma 2, $\mathcal{F}^0$ is an algebraic lower bound for all possible solution curves meeting $\Pi_2$ in angle $\gamma_2$. In this sense, the result of Theorem 4a cannot be improved.

In the other direction we find by a similar discussion:

**Theorem 4*b*.** *Let $\psi_{1,j}^+$ be a decreasing sequence, satisfying (4) and approaching $\psi_1^0$ (see Figure 1a). Then the corresponding normalized forces $\mathcal{F}_j$ satisfy $\mathcal{F}_j = -(1 - \cos \psi_{1,j}^+)$, and $\mathcal{F}_j \nearrow \mathcal{F}^1 = -(1 - \cos \psi_1^0) < 0$ as $j \to \infty$. The value $\mathcal{F}_j$ is achieved at a coordinate $\xi_j$, strictly between the initial plates, with $\xi_j \nearrow \xi_1 =$ the position of $\Pi_1$. The upper bound $\mathcal{F}^1$ is not achieved during this procedure.*

Theorem 4*b* reflects the special role of the barrier **II**, as the least upper height bound of solutions of (2*) in the given family with fixed $\gamma_2$ in $0 < \gamma_2 < \pi/2$, such that the net algebraic force between the plates initially decreases as the plates are brought together with unchanged contact angles. Notably if the starting solution for the procedure were to lie strictly between **II** and **I**, then that force would increase with every decrease in B, cf the discussion in [2] Sec. 6. The solution would rise throughout its traverse until it identifies with **I**, which yields vanishing force; thereafter the force becomes attracting and rises as $O(1/B)$ toward $+\infty$.

We consider finally what happens to lower class neighbors (blue curve in Figure 3). Here the permissible range is (6), subject to the extended definition introduced just below that relation. A reasoning analogous to the one just given shows the existence of a point $\xi_1^-$ interior to the interval, such that the slope of **IV$_0$** at $\xi_1^-$ is exactly $\psi_1^-$. Thus the solution curve becomes a subarc of **IV$_0$** abutting in the intersection point of $\Pi_2$ with the $\xi$–axis, which it cuts in angle $\psi_2 = \pi/2 - \gamma_2$. We are led to:

**Theorem 5.** *As the plate $\Pi_1$ is moved toward $\Pi_2$ with fixed angles $\psi_1^-$ and $\psi_2$, subject to (6), there is a uniquely determined intermediate separation, at which the normalized force $\mathcal{F}^-$ between the plates will be exactly*

$$\mathcal{F}^- = -2(1 - \sin \gamma_2). \tag{10}$$

**Corollary:** *In the course of moving $\Pi_1$ to $\Pi_2$, a neighbor that is initially lower class always achieves exactly at some intermediate position the greatest lower bound of all conceivable algebraic repelling forces obtainable with the given datum $\gamma_2$ on $\Pi_2$.*

The respective positions at which the extreme values are achieved are evaluated in Appendix 5. As we have already pointed out, the lower bound (10) is not achieved by the symmetric solution, for any plate separation.

Note that *we now achieve always the identically same force value, independent of the initial datum $\psi_1^-$ within the admissible range* (6). This value is exactly the upper bound in magnitude that appears for symmetric configurations, and thus *the maximum conceivable repelling force is achieved by every solution in the given neighboring range at a non-zero plate separation, although it is not achieved by the symmetric solution itself.*

The theorems above are supported by computer calculations, in Figure 5. In all upper class cases, as the plates come together the force stays at first repelling, and decreases to exactly the value obtained by inserting the initially imposed value for the inclination ψ into the general force formula (8) for repelling plates. That value approaches the limiting value $-2(1 - \sin\gamma_2)$ for the symmetric case as the initial approximation to the symmetric solution becomes exact. As with the symmetric solution, that limiting value is not actually achieved; instead this solution abruptly (albeit continuously) departs from symmetric behavior as the plate $\Pi_2$ is approached, and the repelling behavior changes rapidly to attracting, with the force becoming positively infinite as $O(1/B)$. The dashed curve gives the forces achieved by the symmetric solution $\mathcal{W}$ at the indicated separations.

The lower class behavior seems to the present authors especially remarkable. *Regardless of the initial choice for ψ within the permissible range, the repelling force magnitude increases to the absolute upper bound $2(1 - \sin\gamma_2)$ for the symmetric case, and in fact achieves that maximum magnitude at a positive plate separation, strictly interior to the range of the motion.* At that crucial point, the behavior begins to reverse, and the forces rapidly become attracting and rise to positive infinity, as in the upper class configuration.

As the initial approximation to $\mathcal{W}$ becomes more precise, the indicated rapid changes all occur within arbitrarily small neighborhoods of $\Pi_2$. The crucial minimizing points are determined to arbitrarily desired accuracy in Appendix 4.

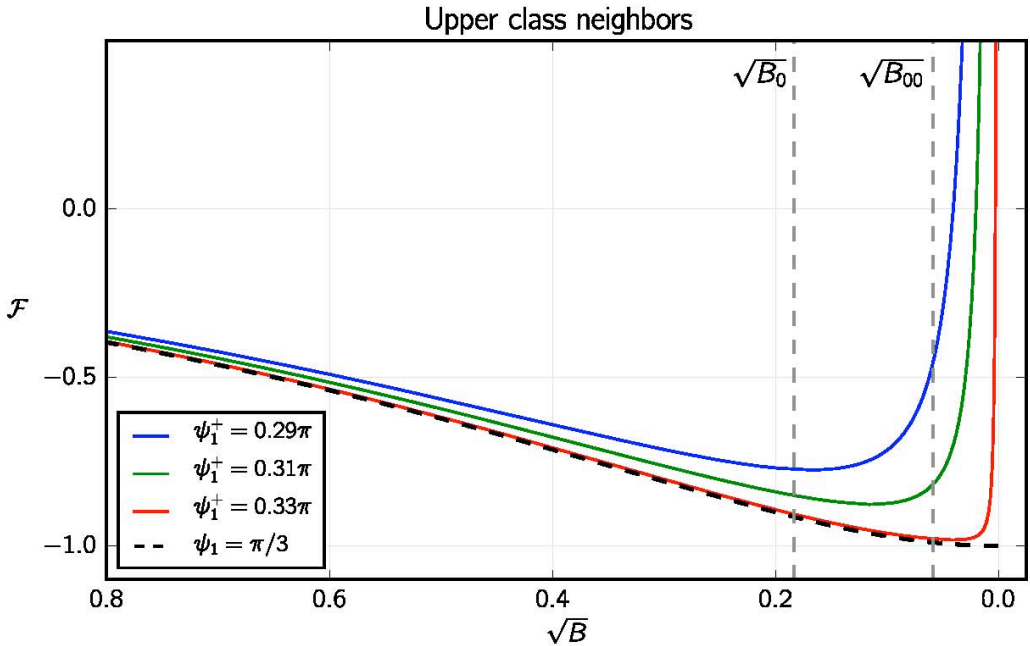

**Figure 5***a*. Forces on the plates as the separation $2B \to 0$, with fixed contact angles, $\psi_2 = \pi/3$; upper class cases. Dashed curve is for symmetric case, with supplementary angles.

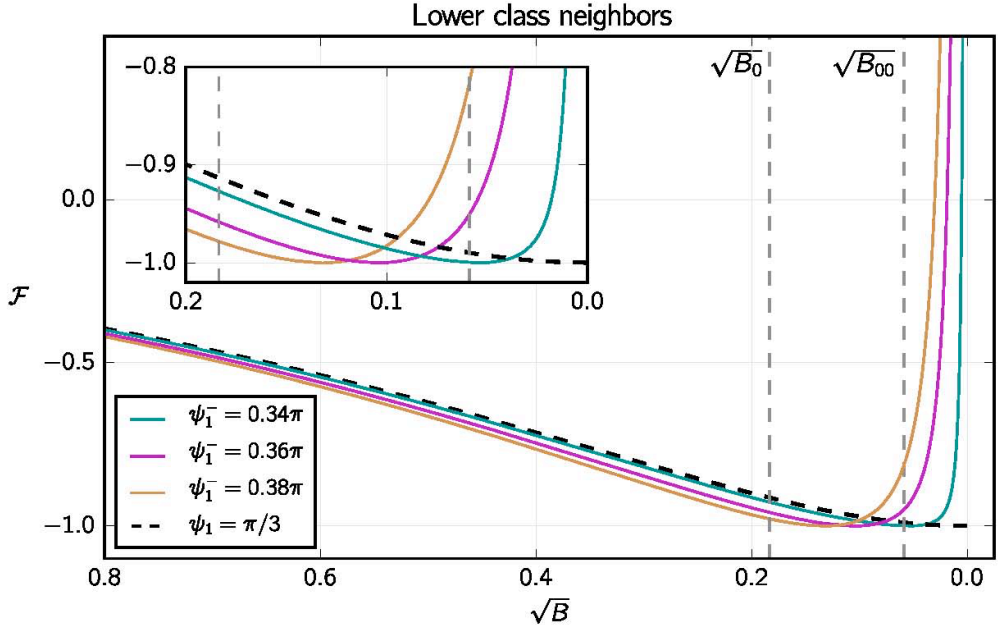

**Figure 5b.** Forces as 2B → 0. with fixed contact angles, $\psi_2 = \pi/3$; lower class cases. Dashed curve is for symmetric case. Note that the minimizing forces are achieved interior to the motion, are identically the same in all cases, and provide a sharp lower bound that is not achieved by any symmetric solution.

**Meniscus heights:** These must in general be obtained by detailed study of the underlying equations (2*a,b), for prescribed angles $\gamma_1$, $\gamma_2$. There is a special interest in asymptotic behavior as the separation parameter B tends to 0 or to ∞. All estimates are contained implicitly within those presented in [1,2,3]; in the interest of conceptual unity, we offer here some explicit discussion for cases of particular interest.

**Case 1.** We consider first a symmetric solution **III**: $\mathcal{W}(\xi;B,\gamma_2)$ (see Figure 2), meeting $\Pi_2$ in angle $\gamma_2$, $0 \le \gamma_2 < \pi/2$, and crossing the $\xi$– axis at $\xi = 0$, the midpoint between the plates. In this case, $\gamma_1 + \gamma_2 = \pi$. We utilize the "first integral"

$$\frac{B}{2}U^2 + \cos\psi = const \tag{11}$$

on any solution curve, from which

$$U^2 = \frac{2}{B}(\cos\psi_0 - \cos\psi) \tag{12}$$

where $\psi_0$ is the angle of crossing with the $\xi$ – axis. From the monotonicity of $\psi$ in $\xi$ we conclude

$$0 < U < \sqrt{\frac{2}{B}(\cos\psi_0 - \sin\gamma_2)} \tag{13}$$

throughout the interval of definition $0 < \xi < 1$. From $(2^*a)$ we find $BU\, d\xi = \cos\psi\, d\psi$, which leads using (12) to

$$\sqrt{2B} = \int_{\psi_0}^{\psi_2} \frac{\cos\psi\, d\psi}{\sqrt{\cos\psi_0 - \cos\psi}} > \int_{\psi_0}^{\psi_2} \frac{\cos\psi\, d\psi}{\sqrt{1-\cos\psi}} = \int_{\psi_0}^{\psi_2} \cot\psi\sqrt{1+\cos\psi}\, d\psi$$
$$> \int_{\psi_0}^{\psi_2} \cot\psi\, d\psi = \ln\frac{\sin\psi_2}{\sin\psi_0} \quad (14)$$

where $\psi_2 = (\pi/2) - \gamma_2$. In view of the monotonicity of $\psi(\xi)$ we conclude immediately that *in the symmetric case, the inclination $\psi$ of the meniscus tends uniformly to $\psi_2$ as $B \to 0$, and thus that asymptotically as the plates come together the meniscus becomes a linear segment at inclination $\psi_0 = \psi_2$.*

**Case 2.** $\gamma_1 + \gamma_2 < \pi$. This is no longer symmetric in the given interval, and is equivalent to $\psi_2 > \psi_1$. Thus by the capillarity equation $(\sin\psi)_\xi = BU$ we find

$$B\int_{-1}^{1} U\, d\xi = \sin\psi_2 - \sin\psi_1 = \lambda > 0 \quad (15)$$

We conclude from (15) the existence of a point $x_m$ within the integration interval, at which $U_m = U(x_m) > 0$. Since $BU^2 + 2\cos\psi$ is constant on the profile curve, we have

$$BU^2 = BU_m^2 + 2(\cos\psi_m - \cos\psi) \quad (16)$$

so that

$$\lambda^2 < 2B^2 \int_{-1}^{1} \left(U_m^2 + \frac{2}{B}(\cos\psi_m - \cos\psi)\right) d\xi$$
$$< 4B^2 U_m^2 + 8B(1 - \sin\gamma_2) \quad (17)$$

from which

$$U_m^2 > \frac{\lambda^2}{4B^2} - \frac{2}{B}(1 - \sin\gamma_2) \quad (18)$$

which becomes positive infinite as the separation parameter B decreases to zero. This tells us in particular that *the curve $U(\xi)$ does not cross the $\xi$ – axis when B is small enough (depending on $\lambda$ and on $\gamma_2$), and thus by (15) there must hold $U(\xi) > 0$ throughout the channel.* Setting now $U_m = \min U$, $U_M = \max U$ in the channel, we obtain from (16)

$$B(U_M^2 - U_m^2) = 2(\cos\psi_m - \sin\gamma_2) < 2(1 - \sin\gamma_2) \quad (19)$$

from which

$$U_M - U_m < \frac{2(1-\sin\gamma_2)}{B(U_M + U_m)} < 2/BU_m \quad (20)$$

when $B < \lambda^2/8(1-\sin\gamma_2)$. In terms of physical heights $u$, (20) asserts that

$$\sqrt{\kappa}(u_M - u_m) < 2/\sqrt{\kappa}u_m, \quad (21)$$

and from (18) we find that

$$\kappa u_m^2 > \frac{1}{\kappa a^2}\left(\frac{\lambda^2}{4} - 4\kappa a^2\right) = \frac{1}{\kappa a^2}\left(\frac{\lambda^2}{4} - 4\text{B}\right). \tag{22}$$

Thus, if we choose B < $\lambda^2/32$, we will have

$$\kappa u_m^2 > \lambda^2 / 8\kappa a^2 \tag{23}$$

and in view of (21) we obtain that *the height change on the traverse of the profile curve tends to zero with the separation, to the same order as the separation.*

**Case 3.** $\gamma_1 + \gamma_2 > \pi$. The discussion is identical to the one just above, except that in this case we find that for small enough B, the meniscus will lie below the $\xi$– axis. The identical magnitude estimates prevail.

In Case 2, there will be a particular value of B at which the meniscus curve coincides with **I**. This is a critical configuration, at which the net horizontal force changes from repelling to attracting, with decreasing B. An analogous configuration appears in Case 3, with **I** replaced by **V**, provided that B < $B_{00}$. If B > $B_{00}$, there are no attracting configurations in Case 3.

With regard to behavior as B → ∞, we observe simply that the height on $\Pi_2$ is obtained explicitly in terms of $\psi_0$ from (12), and that $\psi_0$ can be estimated from the initial integral in (14). By comparison with the solution **I** of Figure 2, one sees that in terms of the rise heights achieved, the plates become rapidly incognizant of each other with increasing separation. Some detailed estimates are given in [2]. The forces become correspondingly vanishingly small, however the transition in direction of force as the datum on $\Pi_1$ traverses that induced by **I** (for which the net horizontal force vanishes) may well be confirmable at large separations.

**Force estimates for narrow channels.** We have shown just above that with the single exception of the symmetric solution **III**, every meniscus meeting $\Pi_2$ in angle $\gamma_2$ with 0 < $\gamma_2$ < $\pi/2$ becomes uniformly ± ∞ in the channel as the separation B → 0, according to the sign of $\gamma_1 + \gamma_2 - \pi/2$. Every such meniscus is attracting when it lies above the barrier **I** or below **V**, and by Lemma 1 above, for any $\gamma_2$ in the range 0 < $\gamma_2$ < $\pi/2$ it suffices to establish that inequality at any single point of the interval. At the right hand end point (on $\Pi_2$), we know the height of **I** to be

$$U_M^{\mathbf{I}} = \sqrt{\frac{2}{\text{B}}(1 - \sin\gamma_2)} \tag{24}$$

We rewrite (16) in the form

$$\text{B}U^2 = \text{B}U_M^2 + 2(\cos\psi_M - \cos\psi) \tag{25}$$

corresponding to the maximum height $U_M$ achieved on $\Pi_2$. Since by monotonicity there holds $\cos\psi > \cos\psi_M$ at all points on the interval, (17) can be replaced by

$$\lambda^2 < 2B^2 \int_{-1}^{1} U_M^2 \, d\xi = 4B^2 U_M^2 \tag{26}$$

and thus to achieve the height $U_M \geq U_M^I$ on $\Pi_2$ it suffices to have

$$\frac{\lambda^2}{4B} \geq 2(1 - \sin\gamma_2). \tag{27}$$

*For all separations B small enough for* (27) *to hold, the (extended) meniscus will lie above* **I** *and will thus yield attracting configurations.* To evaluate the resulting force, we use the relation (1.8) of [3]

$$\mathcal{F} = BU_0^2 \tag{28}$$

where $\mathcal{F}$ is the actual force density normalized by surface tension, and $U_0$ is minimum height of the extended meniscus curve. Using again the first integral $BU^2 + 2\cos\psi = C$, we obtain

$$\mathcal{F} = BU_2^2 - 2(1 - \sin\gamma_2). \tag{29a}$$

In terms of physical parameters, the actual horizontal force per unit length on the contact line orthogonal to the section of the figure is

$$F = \sigma\kappa u_2^2 - 2\sigma(1 - \sin\gamma_2). \tag{29b}$$

For the barrier **I**, the right sides of (29a) and of (29b) vanish. The significance of these relations is emphasized by (21) and (23), which show that the rise heights become asymptotically constant across the channel as the plates come together, and that this constant becomes infinite. Thus from an asymptotic viewpoint, the value $u_2$ in (29b) could be replaced by any of the heights achieved between the plates.

The case $\lambda < 0$ yields to a similar discussion.

**Appendix 1. Determination of** $B_0$: The separation $2B_0$ appears when the plate $\Pi_1$ is tangent to the (barrier) solution **IV₀** at its vertical point; see Figure 2a. We rewrite (2*a,b) as a first order system:

$$\begin{aligned} a) \; & \frac{d\xi}{d\psi} = \frac{\cos\psi}{BU} \\ b) \; & \frac{dU}{d\psi} = \frac{\sin\psi}{BU} \end{aligned} \tag{30}$$

The latter relation separates and can be integrated explicitly. Using the defining property that **IV₀** intersect $\Pi_2$ on the $\xi$–axis at angle $\gamma_2 = \pi/2 - \psi_2$, we obtain when $B = B_0$ that

$$U = \sqrt{\frac{2}{B}(\sin\gamma_2 - \cos\psi)}. \tag{31}$$

At the intersection with $\Pi_1$, we have $\psi = \pi/2$. We place (12) into (11a), and integrate over the known $\psi$– interval. Since the plate separation in the $\xi$- coordinate is exactly 2, we are led to

$$\sqrt{2B_0} = \frac{1}{2}\int_{(\pi/2)-\gamma_2}^{\pi/2}\frac{\cos\psi\, d\psi}{\sqrt{\sin\gamma_2 - \cos\psi}} \tag{32}$$

which determines $B_0$ explicitly. See Figure 6 below.

**Appendix 2. Determination of** $B_{00}$. Here we require the plate $\Pi_1$ to be tangent to **V** at its vertical. In that configuration, **V** coincides with **IV**, and thus **IV** extends as a (negative) solution to $\xi = +\infty$. with $\psi$ tending to 0. Denoting by $U_1$ the normalized intercept height on $\Pi_1$, we find from (11b) that

$$\frac{B_{00}}{2}\left(U^2 - U_1^2\right) = -\cos\psi \tag{33}$$

Allowing the running coordinate $\xi$ to tend to positive infinity, we find $U, \psi \to 0$. Thus

$$U_1 = -\sqrt{2/B_{00}}, \tag{34}$$

see Figure 2b. Using (14) and (15) and observing that the plate separation in normalized coordinates is identically 2, we obtain from (11a)

$$2\sqrt{B_{00}} = \int_{(\pi/2)-\gamma_2}^{\pi/2}\frac{\cos\psi\, d\psi}{\sqrt{B_{00}U_1^2 - 2\cos\psi}}. \tag{35}$$

Placing (34) into (35) yields the explicit formula

$$2\sqrt{B_{00}} = \frac{1}{2}\int_{(\pi/2)-\gamma_2}^{\pi/2}\frac{\cos\psi\, d\psi}{\sqrt{1-\cos\psi}}. \tag{36}$$

The substitution $\psi = 2\tau$ leads to detailed evaluation of (36) in terms of elementary trigonometric functions and logarithms.

Figure 6 displays the critical values $B_0$ and $B_{00}$ in terms of $\gamma_2 = \pi/2 - \psi_2$.

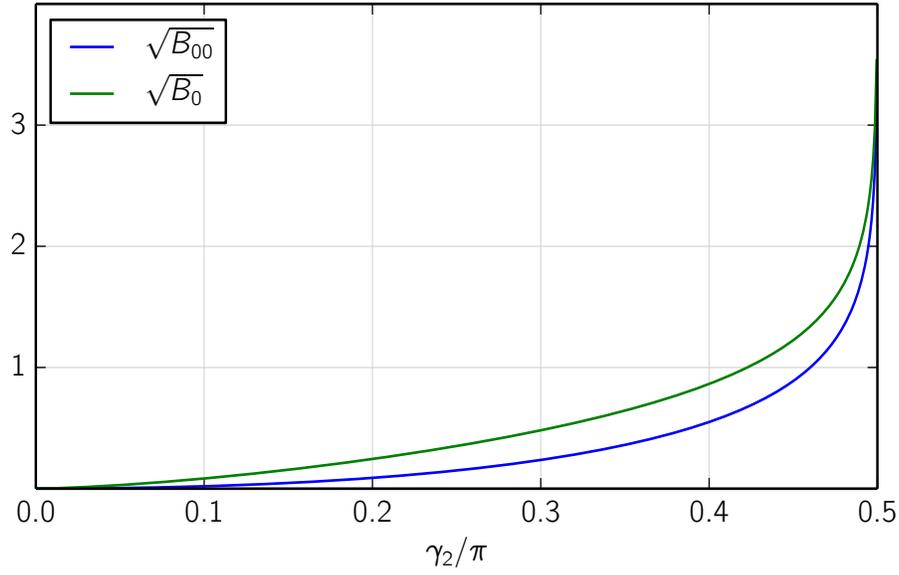

**Figure 6. The critical values $B_0$, $B_{00}$.**

## Appendix 3. Determination of $\psi_1^0$.

See Figure 2a,b. From (11b) we find in the interval between the plates

$$\sqrt{B}U = \sqrt{2\left(\cos\psi_1^0 - \cos\psi\right)}. \tag{37}$$

Since the (non-dimensional) distance between the plates is 2, we may place (37) into (11a) to obtain

$$2\sqrt{2B} = \int_{\psi_1^0}^{\psi_2} \frac{\cos\psi\,d\psi}{\sqrt{\cos\psi_1^0 - \cos\psi}} \tag{38}$$

From (38) we may obtain $\psi_1^0$. The relation is illustrated for four choices of $\psi_2$ in Figure 7 below.

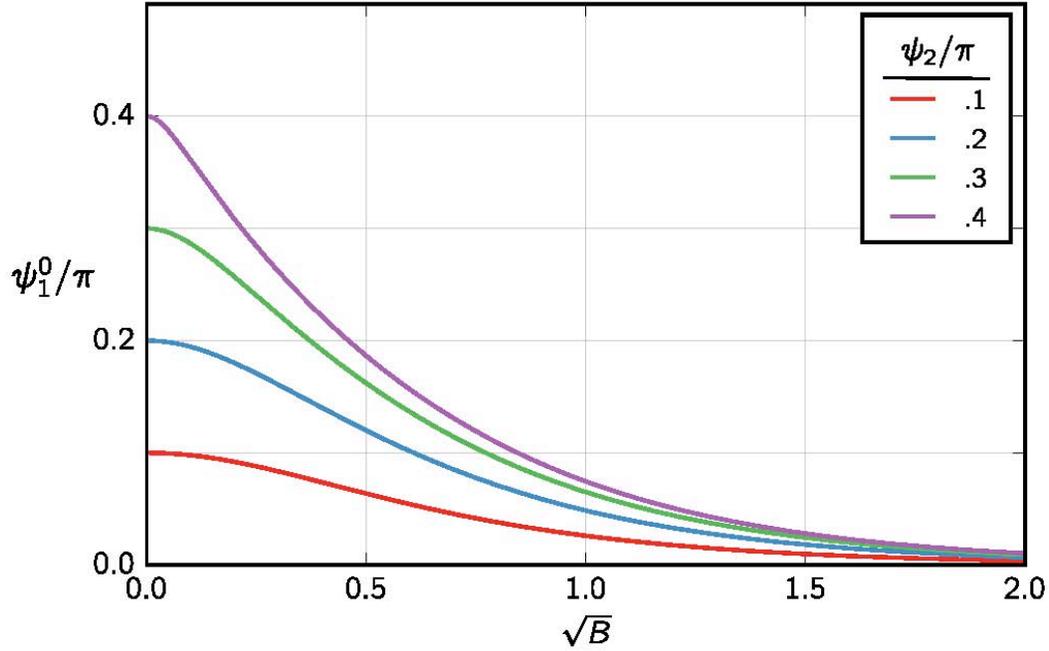

**Figure 7. Determination of $\psi_1^0$.**

**Appendix 4. Determination of $\gamma_{10}$.** This is the "contact angle" in which the solution in the considered family that passes through the intersection of $\Pi_2$ with the $\xi$–axis meets $\Pi_1$. It is defined only when $B \le B_0$, see Appendix 1 where $B_0$ is determined. Thus we assume $B \le B_0$, and we determine the intersection (contact) angle $\psi_{10} = \gamma_{10} - \pi/2$, see Figure 2$b$. From (11$b$) we find that

$$U = \sqrt{\frac{2}{B}} \sqrt{\sin \gamma_2 - \cos \psi} \tag{39}$$

along **IV$_0$**, so that as consequence of the total width being 2 we now obtain using (11$a$) that

$$2 = \frac{1}{\sqrt{2B}} \int_{\psi_2}^{\psi_{10}} \frac{\cos \psi \, d\psi}{\sqrt{\sin \gamma_2 - \cos \psi}} \tag{40}$$

From (40) we may compute $\psi_{10}$. The relation is displayed in Figure 8 below:

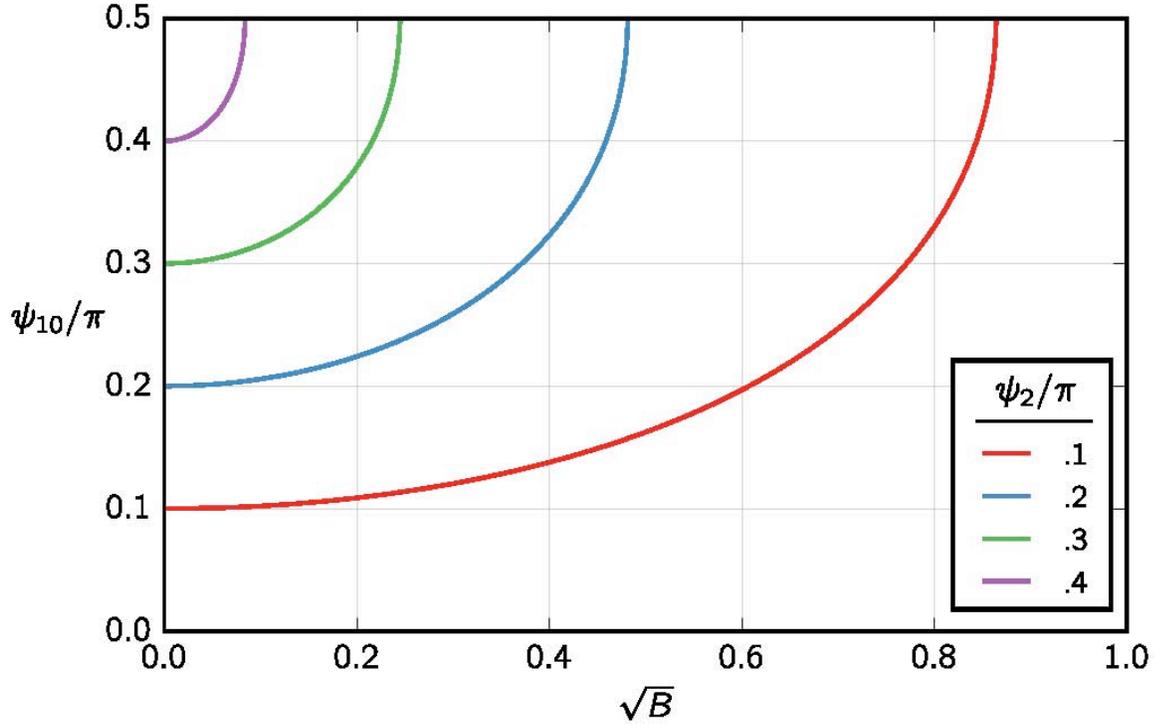

**Figure 8. Determination of $\psi_{10}$.**

**Appendix 5. Extremal positions for repelling configurations.** We have shown both for upper and for lower classes, that for any starting position in the appropriate ranges, there is a unique plate separation, at which the relevant extremal value for the repelling force is actually achieved. It remains to determine the actual values for these separations. Since the reference value of B changes continuously as $\Pi_1$ moves toward $\Pi_2$, we start the discussion by reverting to the original equations prior to normalization:

$$a) \ \frac{dx}{d\psi} = \frac{\cos\psi}{\kappa u} \qquad b) \ \frac{du}{d\psi} = \frac{\sin\psi}{\kappa u} . \qquad (41)$$

**Case 1: Upper class.** Using Theorem 2 and Lemma 2 above, we see that for given $\psi_1^+$ satisfying (4), the extremal configuration is determined as the (unique) solution of (40) cutting the $x$ – axis in the angle $\psi_1^+$ at a point $x_m$, and cutting the $u$ – axis at inclination $\psi_2 =$ $(\pi/2) - \gamma_2$. Using (41b) we obtain that

$$\kappa u^2 = 2\left(\cos\psi_1^+ - \cos\psi\right) \qquad (42)$$

and thus from (41a) the change $\delta x$ from the position $x_m$ to that of $\Pi_2$ becomes

$$\delta x = \frac{1}{\sqrt{2\kappa}} \int_{\psi_1^+}^{\psi_2} \frac{\cos\psi \, d\psi}{\sqrt{\cos\psi_1^+ - \cos\psi}}. \tag{43}$$

Denoting by $2a$ the plate separation at the start of the procedure, we have

$$\frac{\delta\xi}{2} \equiv \frac{\delta x}{2a} = \frac{1}{2\sqrt{2B}} \int_{\psi_1^+}^{\psi_2} \frac{\cos\psi \, d\psi}{\sqrt{\cos\psi_1^+ - \cos\psi}} \tag{44}$$

with $B = \kappa a^2$. For the entire segment, extending from $\Pi_1$ to $\Pi_2$ we find

$$2 = \frac{1}{\sqrt{2B}} \int_{\psi_1^0}^{\psi_2} \frac{\cos\psi \, d\psi}{\sqrt{\cos\psi_1^0 - \cos\psi}}. \tag{45}$$

From (44) and (45) we may eliminate $B$, thus obtaining an explicit expression for the relative change of position corresponding to the chosen $\psi_1^+$. Figure 7 shows the relative extremal positions $\delta\xi$ as function of $\psi_1^+$, for several choices of $\psi_1^0$.

The value $\psi_1^0$ is essential to our discussion, as (together with the explicitly known $\psi_2$) it determines the permissible range for $\psi_1^+$. *We assert that it is determined uniquely in terms of B from* (42), *for any $\psi_2$ in $0 < \psi_2 \le \pi/2$.* To show that, we note initially that the integral

$$\mathcal{J}(\psi_1^0; \psi_2) \equiv \int_{\psi_1^0}^{\psi_2} \frac{\cos\psi \, d\psi}{\sqrt{\cos\psi_1^0 - \cos\psi}} \tag{46}$$

can be made to achieve any value in $(0, \infty)$ by appropriate choice of $\psi_1^0$ in $(0, \psi_2)$. Thus there is at least one admissible solution $\psi_1^0$ of (44), for any $B > 0$. Setting $s = -\cos\psi$, $s_1 = -\cos\psi_1^0$, $s_2 = -\cos\psi_2$, we find

$$\mathcal{J}(s_1; s_2) = -\int_{s_1}^{s_2} \frac{s \, ds}{\sqrt{s - s_1}\sqrt{1 - s^2}}. \tag{47}$$

To obtain uniqueness, it suffices to show that $\mathcal{J}$ is strictly decreasing in $s_1$. We cannot differentiate under the sign, as that leads to singular expressions. $\mathcal{J}$ is nevertheless differentiable in $s_1$, as we see by observing that an integration by parts changes its expression to

$$\mathcal{J} = -\frac{2s_2}{\sqrt{1-s_2^2}} \sqrt{s_2 - s_1} + 2\int_{s_1}^{s_2} \frac{(s - s_1)^{1/2}}{(1-s^2)^{3/2}} ds \tag{48}$$

so that

$$\frac{\partial \mathcal{J}}{\partial s_1} = \frac{s_2}{\sqrt{1-s_2^2}\sqrt{s_2 - s_1}} - \int_{s_1}^{s_2} \frac{1}{(s-s_1)^{1/2}(1-s^2)^{3/2}} ds \tag{49}$$

and since $s_2 < 0$ and $ds_1/d\psi_1^0 > 0$, the assertion follows. Figure 8 shows $\psi_1^0$ as function of B

for several choices of $\psi_2$.

**Case 2: Lower class**

Given an inclination $\psi_1^-$ subject to (6), we must determine the unique point on **IV$_0$** at which this inclination is achieved. As with (21) above, we obtain

$$\frac{\delta\xi}{2} = \frac{\delta x}{2a} = \frac{1}{2\sqrt{2B}} \int_{\psi_2}^{\psi_1^-} \frac{\cos\psi\, d\psi}{\sqrt{\cos\psi_2 - \cos\psi}} \quad (50)$$

for points with inclination $\psi_1^-$ between $\psi_2$ and $\psi_{10} = \gamma_{10} - \pi/2$, and

$$1 = \frac{1}{2\sqrt{2B}} \int_{\psi_2}^{\psi_{10}} \frac{\cos\psi\, d\psi}{\sqrt{\cos\psi_2 - \cos\psi}} \quad (51)$$

when $\delta x = 2a$. (51) determines $\psi_{10}$ in terms of B, and (50) determines the position $\delta\xi/2$ of the point of minimal algebraic force relative to the segment joining the $\Pi_2$ to $\Pi_1$. These relationships are illustrated in Figure 9.

## **Further thoughts**:

(i) We have written little here on attracting solutions. It is not difficult to show that if a solution surface is attracting and if the separation 2*a* is decreased without change in "inner" contact angles, then all new surfaces thus obtained will again be attracting. The converse is true if and only if the menisci are "like". Laplace [10] showed in the case $\gamma_1 = \gamma_2$ that as the plates come together, the attracting force becomes infinite as O(1/B), and that the fluid height rises (or falls) as O(1/√B). More precise estimates for a general case appear in [2] and in [3]. Expressions for the forces relevant to partly immersed bodies of general shape can be found, e.g., in Miersemann [8].

(ii) The above material on meniscus heights relates to a paper of Vella and Mahadevan, see the discussion in [11].

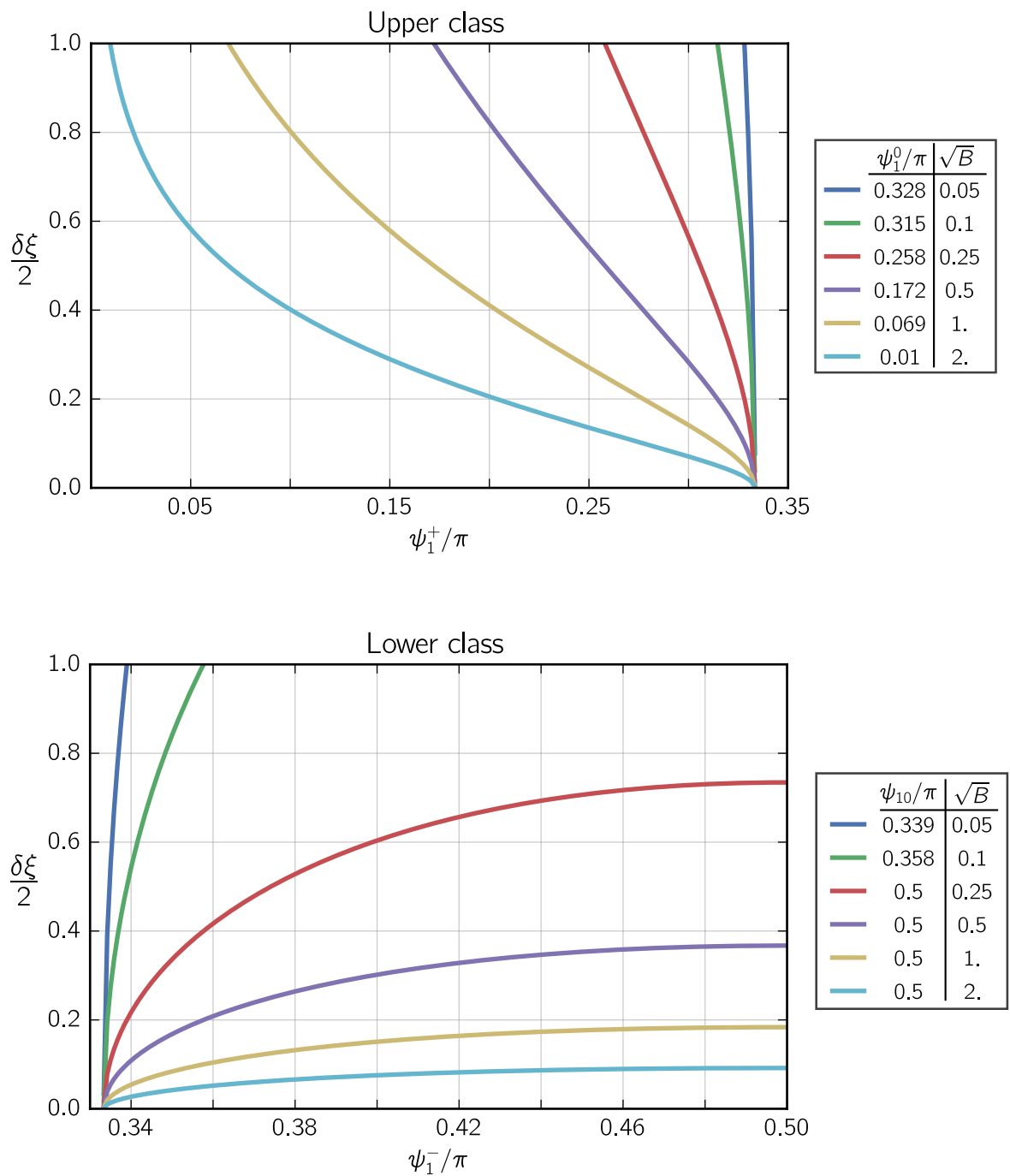

**Figure 9.** Determination of position of extremal points, relative to $\Pi_2$; upper and lower classes, $\psi_2 = \pi/3$.

## Acknowledgments.

The second author is indebted to the Mathematische Abteilung der Universität, in Leipzig, for its hospitality and for excellent working conditions, during his visit of October 2014. He is indebted to the Simons Foundation for a very helpful grant facilitating his travel. He thanks Erich Miersemann for incisive comments on some conceptual points related to the problem considered.

Department of Biochemistry and Biophysics,
California Institute for Quantitative Biosciences,
University of California
San Francisco, CA 94143, USA
raj.bhatnagar@gmail.com

Mathematics Department
Stanford University
Stanford, CA 94025, USA
finn@math.stanford.edu

# Addenda to the Preceding Paper

Rajat Bhatnagar and Robert Finn

**Summary:** This work contains largely afterthoughts, relating to the paper immediately preceding it. We correlate and interpret our contributions in that paper, relative to those of an earlier publication by Aspley, He and McCuan. We propose specific laboratory experiments, suggested by formal predictions of those two papers.

The material of this note can be viewed as a continuation of the study in our immediately preceding paper; accordingly we adopt here generally the same notations, equation numerations, figure designations and citation list. Designations for items in those categories that are new to the present work will be prefixed by the symbol "*A*". The directly preceding work itself is designated here as the item [*A*-1] in our reference list.

## Addendum 1

We begin by interpreting the continuing contributions of the present authors in the context of the recent contribution [4].

A direct comparison of the contributions strikes us as something of a challenge: In both approaches the contact angles $\gamma_1$, $\gamma_2$ on the inward facing plate sides are prescribed, as is the plate separation. Consequent information is sought, on the form of the free interface and on the net forces that act on the plates. In the present authors' approach, the focus is on the profile contours, of which the changes are observed for arbitrary but fixed choice of $\gamma_2$, as $\gamma_1$ traverses all admissible values, for varying choices of separation. Forces are treated as subsidiary information, which can be determined in terms of information established for the profile. Solutions are organized into categories, each demarcated relative to adjacent ones with differing geometrical properties, in terms of reference "barrier" solutions, that serve as guideposts to the varying sorts of behavior.

In the approach of [4], the emphasis is on the forces, which the authors calculate directly in terms of the plate separation, for three different categories of contact angle choices. The interface profiles themselves become subsidiary, and are not examined in explicit detail of their dependence on underlying parameters of the data.

Both approaches cover all configurations that can occur; however in our procedure ten geometrically distinct configurational choices appear; in the approach of [4], only three underlying distinctions are evident. Thus there are geometric distinctions appearing that are not directly recognized by the latter procedure. We will point them out below.

In the sense of logical development of the theory, [4] presents in our view a central and basic contribution, in establishing that every configuration can be assigned to one of exactly three kinds of force/separation profiles. It does not clearly delineate however all geometric changes to the free interface that can occur. As a particular instance, we see no way to locate within that theory the critical plate separation $B_{00}$, at which the set of all attracting solution curves between the plates changes from simply to doubly connected. Within the "barrier" structure as developed in the earlier citations [1,2,3], the changing point is prominent in the theory and easily calculated.

It is clearly desirable to correlate the two procedures to the extent that can be done, and we make an initial start in that direction in the material below. The details of what happens can be striking and varied, and we found all ten distinct kinds of behavior surprisingly accessible to detailed description. Each of the eight barriers in Figure 2 delineates a physical demarcation of individual and idiosyncratic interest, relative to the "separation parameter" $B = \kappa a^2$. We list them as follows, assuming a downward gravity field and a representative $\gamma_2$ on $\Pi_2$, in the range $0 \leq \gamma_2 < \pi/2$.

$\mathcal{T}$: **The solution obtained for datum** $\gamma_1 = 0$. See Figure 2. To characterize the solution formally, we are guided by the discussion in [A-1] starting from (11) onward. We find respectively, on the plates $\Pi_1$, $\Pi_2$, the boundary conditions

$$\psi_1 = -\pi/2, \quad \psi_2 = \pi/2 - \gamma_2. \tag{A1}$$

For the equivalent first order system

$$\xi_\psi = \frac{\cos\psi}{BU}, \quad U_\psi = \frac{\sin\psi}{BU} \tag{A2a,b}$$

of the interface equation

$$(\sin\psi)_\xi = BU \tag{A3}$$

we can integrate (A2b) to obtain the "first integral"

$$\frac{B}{2}U^2 + \cos\psi = C \tag{A4}$$

on any solution curve, and thus, using (A4) in (A2a) and integrating over the interval between the plates, we are led to

$$2 = \frac{1}{\sqrt{2B}} \int_{-\pi/2}^{\pi/2-\gamma_2} \frac{\cos\psi}{\sqrt{C - \cos\psi}} d\psi. \tag{A5}$$

Here we have introduced the data (A1). One sees immediately that for given B, (A5) uniquely determines $C$. The solution $\mathcal{T}$ is thus uniquely determined by B and by the data (A1). The solution can be written in the indicated parametric form, by solving (A4) for a (positive) $U(\psi;B)$, which is then inserted into (A2a,b).

$\mathcal{G}$: **The solution obtained for datum** $\gamma_1 = \pi/2$. The discussion is identical to that for $\mathcal{T}$, however the datum for $\psi_1$ in (A1) now becomes $\psi_1 = 0$.

**I:** **The (unique) solution meeting** $\Pi_2$ **in angle** $\gamma_2$ **as above and extending as a graph to** $x = -\infty$. We obtain this solution by setting $C = 1$ in (A4).

**II:** **The (unique) solution meeting** $\Pi_2$ **in angle** $\gamma_2$ **and passing through the intersection of** $\Pi_1$ **with the** $\xi$ – axis. The requirement on $\Pi_1$ yields $C = \cos\psi_1$ (from (A4)), and $\psi_1$ is then determined by the requirement

$$2 = \frac{1}{\sqrt{2B}} \int_{\psi_1}^{\pi/2-\gamma_2} \frac{\cos\psi}{\sqrt{\cos\psi_1 - \cos\psi}} d\psi \tag{A6}$$

cf the discussion in Appendix 5 of [$A$-1].

**III:** **The (symmetric) solution meeting** $\Pi_1$ **in angle** $\gamma_1 = \pi - \gamma_2$. This solution crosses the $\xi$–axis at the midpoint between the plates, and we must take account that the orientation above that axis is the reverse of the one below. Restricting attention to the right half, in which $U > 0$, and denoting the crossing angle with $\psi_0^*$, we find $C = \cos\psi_0^*$, and

$$1 = \frac{1}{\sqrt{2B}} \int_{\psi_0^*}^{\pi/2-\gamma_2} \frac{\cos\psi}{\sqrt{\cos\psi_0^* - \cos\psi}} d\psi \tag{A7}$$

which determines $\psi_0^*$. The solution is then determined explicitly using (A2*a,b*) as above.

**IV:** **The solution obtained for data $(\pi, \pi/2 - \gamma_2)$ on $(\Pi_1, \Pi_2)$. This solution clearly lies below III.** It is the lower analogue of $\mathcal{T}$, in being the lowest solution joining the plates that can be achieved for the prescribed $\gamma_2$. It will contain points lying directly above corresponding points of the ensuing profiles **IV$_0$** and **V** if and only if those profiles do not extend to meet both plates. In turn, those events depend essentially on the plate separation, see the **Cases I, II, III** below.

**IV$_0$:** **The solution meeting $\Pi_2$ on the $\xi$–axis, at inclination $\psi_2 = \pi/2 - \gamma_2$.** This surface profile is rigidly attached to $\Pi_2$ and independent of plate separation. It extends across the channel to $\Pi_1$ only if $B \leq B_0$, see the designated cases below.

**V:** **The lower counterpart to I:** This profile lies in the lower half-plane, meets $\Pi_2$ in angle $\gamma_2$ and extends to $\xi = +\infty$. It is rigidly attached to $\Pi_2$ and extends to $\Pi_1$ only if $B \leq B_{00}$.

Observed behavior depends in essential ways on the relative positions of **IV,** of **IV$_0$** and of **V.** These change according to the separation B, as indicated in Figure *A*–1 below:

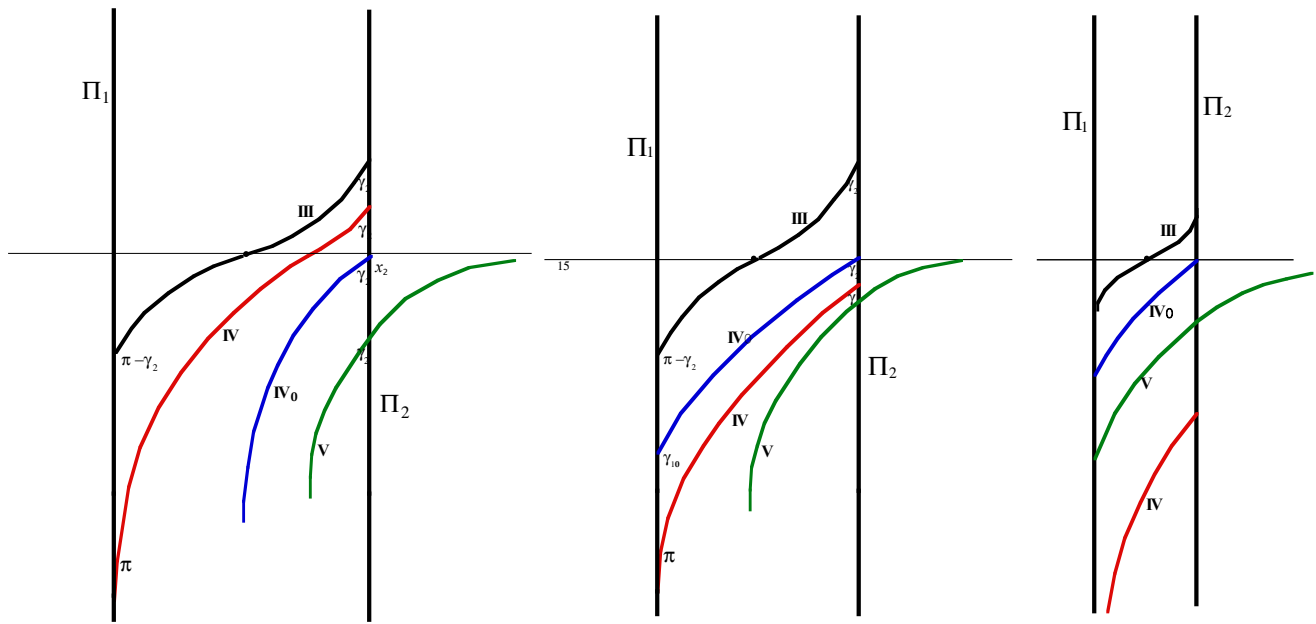

    Case 1. Wide: $B > B_0$      Case 2. Intermediate: $B_0 \geq B > B_{00}$      Case 3. Narrow: $B_{00} \geq B$

**Figure *A*–1. Configurations adopted by characteristic interfaces below III, for varying plate separations (not to scale)**

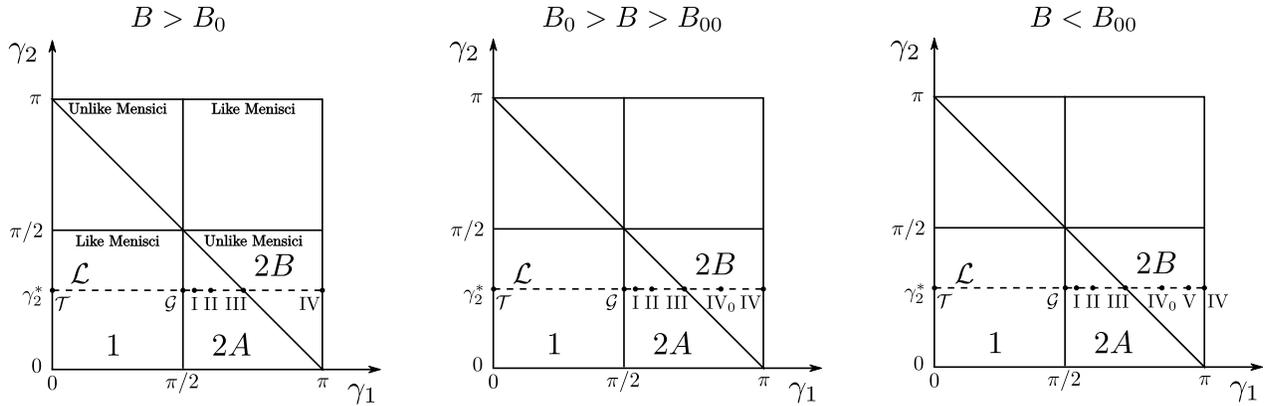

**Figure $A$ –2.** Domain choices for data pairs. The three figures together provide an analogue of Figure 3 of [4], modified for the present discussion; $\gamma_1$ and $\gamma_2$ are interchanged. The dashed line $\mathcal{L}$ is the restriction to fixed $\gamma_2^*$ on $\Pi_2$. The dots on $\mathcal{L}$ represent corresponding barriers in Figure 2.

The above barrier solutions relate to the framework introduced in Theorem 4 of [4], see Figure $A$ -2 above. That figure is the analogue in the present notation (continued from earlier work by the present authors), of the Figure 3 in [4]. Some of the interpretations follow:

We note initially, that the points in the closed lower halves of the large squares in Figure $A$ -2 correspond bi-uniquely to the set of all distinct solution curves of (2*a,b*) that join the two plates (Essentially Sec. II of [1]). The one –parameter subset of those curves arising from the initial choice for $\gamma_2$ is indicated by the line $\mathcal{L}$: $\gamma_2 = \gamma_2^*$ of Figure $A$ -2.

All barrier solutions correspond to individual points on $\mathcal{L}$. Using notation introduced just following Figure 2, we find the following interpretations, for general configurations:

$\mathcal{T}$: *no point of any solution arising from $\mathcal{L}$ can lie above $\mathcal{T}$.*

$\mathcal{R}_{\mathcal{T}\text{-}\mathcal{G}}$: *The set of meniscus profiles lying on or below $\mathcal{T}$ and on or above $\mathcal{G}$. They are exactly the solutions from the part of $\mathcal{L}$ lying in the square "**1**" of Figure $A$– 2.*

**I:** *this is the (unique) positive solution yielding zero force between the plates. In the "wide" configuration there is no negative solution in that category. Negative solutions with zero force do not appear until the plate separation is small enough, so that $B \leq B_{00}$. Note that all menisci from points interior to "**1**" are "like", in the sense that $(\pi/2 – \gamma_1)(\pi/2 – \gamma_2) > 0$.*

$\mathcal{R}_{\mathcal{G}\text{-}\mathbf{I}}$: *the set of solution curves with "unlike menisci" for which the plates attract each other. These can also be characterized as the set of positive attracting solutions, for which the points of minimum height lie exterior to the plates.*

**II:** *this solution separates the repelling solutions that cross the x – axis between the plates, from those crossing the axis to the left of the plates. In the former case the negative minimizing force will be achieved by decreasing the separation; in the latter case one must increase the separation, cf. the discussion in Sec. 7 of [2].*

$\mathcal{R}_{\text{I-II}}$: *the set of repelling solutions crossing the x – axis to the left of the plates. For a fixed $\gamma_2$ in $0 \leq \gamma_2 < \pi/2$, this corresponds bi-uniquely to the indicated horizontal line segment of $\mathcal{L}$ in the region "2A" of Fig. A – 1. From another point of view, this is exactly the set of repelling solutions situated above the symmetric one III, for which one must widen the channel to achieve the algebraic minimum of the repelling forces.*

**III:** *This "symmetric" solution ($\gamma_1 + \gamma_2 = \pi$) is the unique solution that remains repelling as the plates come together with fixed contact angles. With decreasing plate separation, the repelling force increases monotonically in magnitude to the limiting density $\mathcal{F} = -2\sigma(1 - \sin\gamma_1) = -2\sigma(1 - \sin\gamma_2)$.*

$\mathcal{R}_{\text{II-III}}$: *the set of repelling solutions crossing the x – axis between the plates, and lying above III. This set corresponds bi-uniquely to the part of $\mathcal{L}$ that lies between II and III in "2A". This is also exactly the set of "upper-class neighbors" introduced in [A-1]. They are the repelling solutions lying above III, with the property that by narrowing the channel one achieves an algebraic minimum of repelling forces. This minimum is $\mathcal{F} = -2\sigma(1 - \sin\gamma_1)$, where $\gamma_1$ is the chosen initial datum within the indicated range.*

The further material varies according to the separation criteria. See Figure *A*–1 above.

**Case 1:** $B > B_0$. Then **IV** lies above $\mathbf{IV_0}$, which does not extend to $\Pi_1$. We find:

$\mathcal{R}_{\text{III-IV}}$: *the set of all solutions joining the plates, and lying below III. These solutions all repel, and arise from a portion of $\mathcal{L}$ lying within "2B". There are no attracting solutions below III.*

**Case 2:** $B_0 > B > B_{00}$. In this event, $\mathbf{IV_0}$ lies above **IV** and extends to $\Pi_1$, while on $\Pi_2$ **IV** still lies above **V**, which does not connect the plates. We obtain:

$\mathcal{R}_{\text{III-IV}_0}$: *the set of all solutions joining the plates, lying below III, and crossing the x – axis between the plates. This is identically the set of "lower-class neighbors" introduced in [A-1]. They are the repelling solutions lying below III, with the property that by narrowing the channel one achieves an algebraic minimum of repelling forces. This minimum now becomes $\mathcal{F} = -2\sigma(1 - \sin\gamma_2)$, and is the absolute minimum for all conceivable solutions, extending to Both plates and achieving the prescribed $\gamma_2$ on $\Pi_2$. In contrast to upper-class solutions, this value is independent of the choice of starting element in the set.*

$\mathcal{R}_{\text{IV}_0\text{-IV}}$: *the set of all solutions joining the plates, lying below $\mathbf{IV_0}$. Alternatively it is the set of all solutions joining the plates and crossing the x – axis to the right of $\Pi_2$.*

$\mathcal{R}_{\text{III-IV}}$: *the set of all solutions joining the plates and lying below III. Within this set, $\mathbf{IV_0}$ provides the algebraic minimum $\mathcal{F} = -2\sigma(1 - \sin\gamma_2)$ for the force density, as the plates come together with fixed contact angles. All these solutions correspond to a segment on $\mathcal{L}$ within "2B".*

**Case 3:** $B < B_{00}$. Now **V** lies between $\mathbf{IV_0}$ and **IV** and connects the plates, and new attracting solutions appear below **V**. $\mathcal{R}_{\text{III-IV}_0}$ is defined as in Case 2. However new sets of interest appear:

$\mathcal{R}_{\text{IV}_0\text{-V}}$: *the set of repelling solutions joining the plates, and lying below $\mathbf{IV_0}$. Again it is*

*alternatively the set of all solutions joining the plates and crossing the x – axis to the right of $\Pi_2$. Within this set, there is no algebraic minimum for the horizontal force.*

$\mathcal{R}_{V-IV}$: *the set of all solutions below **V** and joining the plates. These provide, for small separations, a second set of solutions with unlike menisci, which nevertheless attract. These solutions correspond to a final segment of the line $\mathcal{L}$ within* "2B"

The varying behavioral configurations are exhibited graphically, in Figure *A*–3 below. Note especially the second "attracting" zone that appears when B < $B_{00}$.

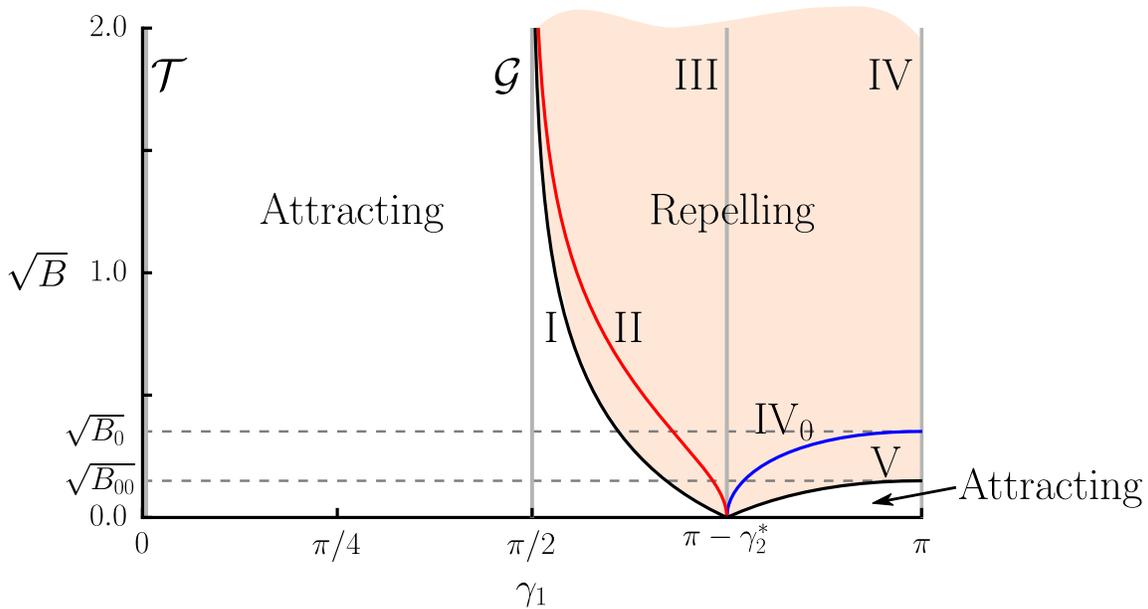

**Figure *A*–3.** Graphical representation of the various configurations achievable by changing plate separation and contact angle $\gamma_1$, for prescribed $\gamma_2 = \pi/4$.

**Comments:** In all cases considered, the presented material is based directly on the original (nonlinear) equations of capillarity, as formulated by Thomas Young in 1805. There are no linearizing approximations.

We note especially the transition that occurs when B decreases past $B_{00}$. When B > $B_{00}$, the set of attracting solutions is simply connected as a point set in the ($\xi$,,$U$) plane, and lies above the $\xi$– axis. There is a discontinuous transition in behavior as B moves below $B_{00}$, in the sense that for *all* B < $B_{00}$, *a second component of the attracting set appears, lying below the horizontal axis. It is disjoint from the original set, and consists entirely of "unlike" menisci, in the sense that the plate contact is wetting on $\Pi_2$, but nonwetting on $\Pi_1$.*

As already remarked, the representation based on the Figure *A*-2 is remarkable in establishing that in an underlying sense, only three kinds of qualitative behavior can occur. Nevertheless, ten distinct cases are distinguished above, each correlating to distinct and observable delineations of physical criteria; that is, the distinctions introduced in [4] subdivide into further cases of individual interest. Further, [4] distinguishes only the qualitative form of the "force profile" as determined by the separation, and does not in itself provide information on the lateral position of the minimizing (or crossing) point. Those lateral positions can however be decisive in several senses, some of which are indicated in the above descriptions.

# Addendum 2

We propose here four specific laboratory experiments, designed to display some of the idiosyncratic behavior indicated above, in concrete physical configurations.

**E1:** As pointed out in [6], the barrier **I**, which is one of the four solutions of (2a,b) meeting $\Pi_2$ in prescribed angle and extending to infinity, itself yields zero force between any two vertical plates meeting it at the incident contact angles, and it separates continua of repelling solutions from continua of attracting solutions.

From (2a,b) we obtain after some manipulation the explicit parametric form for **I**:

$$\kappa u^2 = 2(1 - \cos\psi)$$

$$\sqrt{\kappa}(x - x_2) = \log\frac{\tan(\psi/4)}{\tan(\psi_2/4)} - 4\sin\frac{\psi + \psi_2}{4}\sin\frac{\psi - \psi_2}{4}$$

$$0 < \psi < \psi_2 < \pi/2$$

Keeping $\Pi_2$ fixed and letting $\Pi_1$ move out toward $-\infty$, **I** remains unaffected, and we obtain plate configurations of arbitrarily large separation, for which **I** provides explicit solutions yielding zero force. Small changes in the contact angle $\gamma_1$ then provide corresponding changes from attracting to repelling configurations or vice versa, according to the direction of change. Although the magnitudes of the forces tend rapidly to zero with increasing separation, it should be feasible to measure the change in direction in a controlled experiment. For example, one could leave $\Pi_1$ free to move laterally, and observe the direction of motion.

Figure *A*–4 below provides an indication of what to expect. The curves above **I** are attracting, and prolongations do not cross the *x*–axis. The curves below **I** are repelling. The identical comments apply, for arbitrary separation of the plates, and the horizontal forces are independent of the separation, if the curves are prolongations of the ones indicated. The actual forces in the indicated configuration may be quite small, but it should not be difficult to measure their directions.

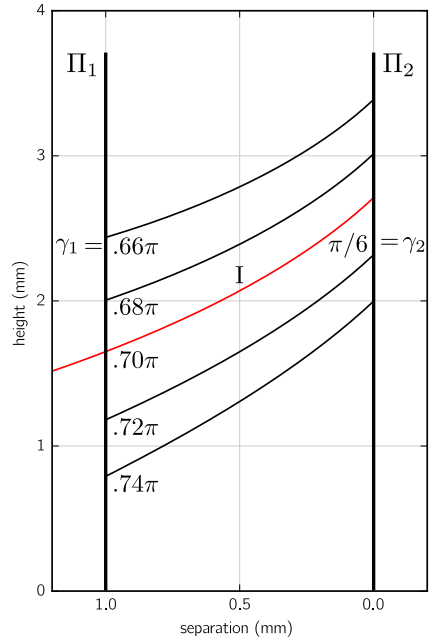

**Figure *A*–4: the zero-force curve I in a particular configuration. The profiles above I yield net attracting force, those below I yield repelling force.**

**E2:** There are multiple opportunities for exhibiting unlike menisci that provide attracting forces. Although the existence of such entities conflicts with generally accepted views deriving from cumulative professional experience, there are in fact more ways for that to occur than there are for appearance of like menisci. Assuming as above that $0 < \gamma_2 < \pi/2$ (otherwise simply reflect in the $x$ – axis) then "like menisci" cannot occur when $u < 0$. However the region $\mathcal{R}_{V\text{-}IV_0}$ becomes a non-null open set simply covered by attracting solutions with unlike menisci, whenever $B < B_{00}$. Additionally, when $u > 0$, the set $\mathcal{R}_{\mathcal{G}-I}$ has the requisite properties for any B.

Since all solutions are (up to routine calculations) explicitly known, it will not be difficult to choose convenient parameters, to establish the relevant behaviors in all configurations under laboratory conditions, and compare directly with prediction.

**E3.** Figure 5*b* (note especially the insert) displays the remarkable property of "lower class neighbors" that for *any* solution in the "lower class" region, the family of solutions obtained by displacing $\Pi_1$ toward (or away from) $\Pi_2$ includes the particular solution **IV₀**, for which the minimum conceivable algebraic repelling force $\mathcal{F} = -2\sigma(1 - \sin\gamma_2)$ corresponding to the datum $\gamma_2$ is strictly achieved. A procedure for determining this critical separation is idescribed in Appendix 5 of A –1. All relevant quantities are explicit, and are subject to direct experimental verification. One would want to verify especially the independence of the calculated minimal algebraic force, relative to the particular choice of solution in the "lower class" range.

**E4.** The following observations strike us as having a special experimental interest. We focus attention on the barriers **IV₀** and **V**, noting that for fixed $\gamma_2$ they are both tied rigidly to $\Pi_2$. We consider the particular solution **IV₀** on a small interval $\mathcal{I}_\delta$: $0 < x_2 - x \leq \delta$ adjacent to $\Pi_2$, see Figure 3. We introduce a vertical plate $\Pi^*$ at the left hand endpoint $x_2 - \delta$, meeting **IV₀** at a point $p$ that tends to ($x_2$, 0) as $\delta \to 0$. At $p$, the "contact angle" of **IV$_0$** with $\Pi^*$ will be $\gamma^*$, approaching $\gamma_2$ as $\delta \to 0$, with $h^* \to 0$.

Since **IV$_0$** and **V** have the same inclinations at $x_2$, and since **V** has more negative curvature than does **IV$_0$** at all points in $\mathcal{I}_\delta$, the angle $\gamma^*$ with the vertical must appear on **V** at an intermediate point ($x^*$, 0) of $\mathcal{I}_\delta$. (This point can be found explicitly, see Appendix 5).

Now fix the angle $\gamma^*$, and shift $\Pi^*$ toward $\Pi_2$. The segment of **IV$_0$** between $\Pi^*$ and $\Pi_2$ will be replaced by another solution in a narrower channel, but with the unchanged boundary angles. When the point ($x^*$, 0) is attained, these angles will coincide with those of **V.** By the general uniqueness theorem, the solution curve has shifted downward a distance exceeding $h$, and the net (repelling) force has changed from $-2\sigma(1 - \sin\gamma_2)$ to zero. We find easily $\sqrt{\kappa}h = \sqrt{2(1 - \sin\gamma_2)}$, which under customary earthbound conditions is a length readily accessible to measurement.

In general terms, the anticipated experimental result is that *for narrow channels, a fixed non-zero change of boundary height and of corresponding repelling force can be attained by arbitrarily small change of boundary angle or of plate separation.* The phenomenon can be correlated to the sharp changes exhibited in the portion to the lower right of Figure 5*b* in A-1.

All quantities described here by terms such as "small" or "narrow" can be made explicit using procedures such as in Appendix 5. It is difficult for us to interpret the significance of the result, especially in view of established difficulties in experimental determination of contact angle. The result indicates that, at least for repelling configurations, small changes in contact angle can result in relatively large changes of height and of force. The results of carefully controlled experiments could have a special interest and lead to new insights.

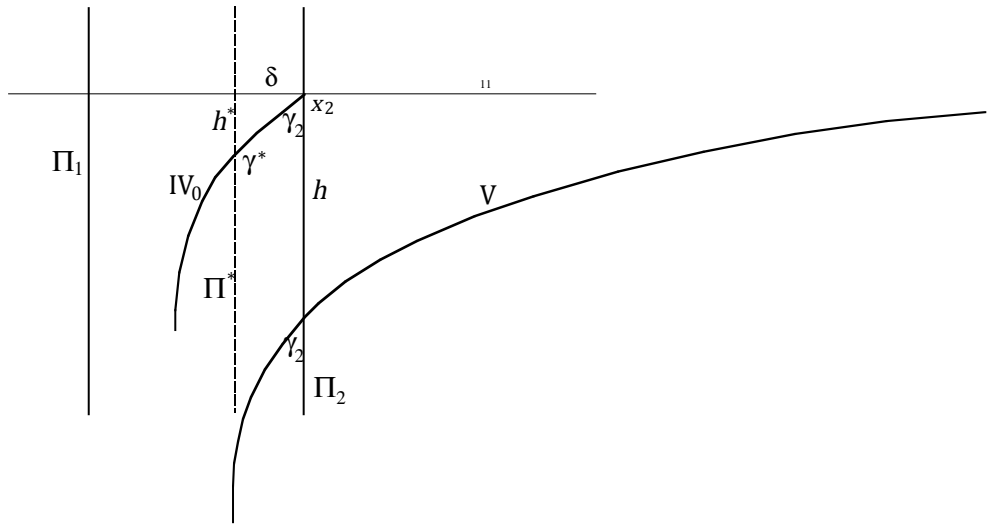

**Figure A–5.** **Abrupt (nearly discontinuous) dependence of height and of repelling force on the contact angle**

**Acknowledgments.** The second author is indebted to the Mathematische Abteilung der Universität, in Leipzig, for its hospitality and for excellent working conditions, during his visit of October 2014. He is indebted to the Simons Foundation for a very helpful grant facilitating his travel. He thanks Erich Miersemann for incisive comments on some conceptual points related to the problem considered.

Department of Biochemistry
        and Biophysics,
California Institute for
        Quantitative Biosciences,
University of California
San Francisco, CA 94143, USA
raj.bhatnagar@gmail.com

Mathematics Department
Stanford University
Stanford, CA 94025, USA
finn@math.stanford.edu